\documentclass[12pt,reqno,twoside,final]{article}
\usepackage{amsmath,amsthm,amsfonts,amssymb,setspace,verbatim,cite,fancyhdr}

\usepackage[english]{babel}
\usepackage[utf8]{inputenc}
\usepackage{times}

\usepackage[T1]{fontenc}
\usepackage{mathrsfs}

\let\oldnsubseteq\nsubseteq 
\usepackage{mathabx} 
\let\nsubseteq\oldnsubseteq

\usepackage{enumitem}
\usepackage{accents}
\usepackage{fge}

\usepackage[top=15mm,bottom=13mm,left=22.8mm,right=22.8mm,includeheadfoot,headsep=8mm,footnotesep=5mm]{geometry} 

\usepackage[hypertexnames=false]{hyperref} 
\hypersetup{
pdfencoding=auto,
pageanchor=true,
plainpages=true,
pdftoolbar=true,
pdfmenubar=true,
pdfpagemode=UseOutlines,%(show bookmarks);UseNone;UseThumbs (show thumbnails);FullScreen;UseOC (PDF % 1.5);UseAttachments(PDF 1.6)
pdffitwindow=false,
pdfstartview={FitH},
pdfnewwindow=true,
colorlinks=true,
linkcolor=black,
citecolor=black,
filecolor=black,
urlcolor=black,
linktoc=all,
breaklinks=true,
}

\usepackage{mdframed,xcolor}
\usepackage{xpatch}

\makeatletter
\xpatchcmd{\endmdframed}
  {\aftergroup\endmdf@trivlist\color@endgroup}
  {\endmdf@trivlist\color@endgroup\@doendpe}
  {}{}
\makeatother

\mdfdefinestyle{defStyle}
{linewidth=0.4mm,
skipabove=0.5\baselineskip,
skipbelow=0.9\baselineskip,
splittopskip=\topskip,
innertopmargin=0pt,
theoremspace={},
}

\newmdtheoremenv[linewidth=0.4mm,skipabove=0.9\baselineskip,skipbelow=0.9\baselineskip,splittopskip=\topskip,%
innertopmargin=0pt]{theorem}{Theorem}[section]
\newmdtheoremenv[linewidth=0.4mm,skipabove=0.5\baselineskip,skipbelow=0.9\baselineskip,splittopskip=\topskip,%
innertopmargin=0pt]{lemma}[theorem]{Lemma}
\newmdtheoremenv[linewidth=0.4mm,skipabove=0.5\baselineskip,skipbelow=0.9\baselineskip,splittopskip=\topskip,%
innertopmargin=0pt]{corollary}[theorem]{Corollary}
\newmdtheoremenv[linewidth=0.4mm,skipabove=0.5\baselineskip,skipbelow=0.9\baselineskip,splittopskip=\topskip,%
innertopmargin=0pt]{proposition}[theorem]{Proposition}
\mdtheorem[style=defStyle]{definition}[theorem]{Definition}
\newmdtheoremenv[style=defStyle]{hypothesis}[theorem]{Hypothesis}

\usepackage{thmtools}

\declaretheorem[style=definition,name=Remark,sibling=theorem]{remark}

\numberwithin{equation}{section}

\makeatletter
\def\th@plain{%
  \thm@notefont{\bfseries}% same as heading font
  \itshape % body font
}

\def\th@definition{%
  \normalfont % body font
  \thm@notefont{\bfseries}% same as heading font
}
\makeatother

\renewcommand{\thefootnote}{\fnsymbol{footnote}}

\newcommand{\makepapertitle}{
\pdfbookmark[1]{Title page}{Title_page}
\thispagestyle{empty}\vspace*{8mm}
\begin{spacing}{1.1}
\LARGE\MakeUppercase{\papertitle}
\end{spacing}\vspace*{10mm}
}

\newcommand{\paperfirstauthor}{%
{\large\sc\firstauthor}\\[2mm] %\FootCorAuth%\footnote{\firstthanks}} \\[2mm]
{\rm\firstaddress \\ 
\firstemail}\\[8mm]
}

\newcommand{\makeabstract}{%
\begin{minipage}{0.9\textwidth}
{\small {\sc Abstract.}
 \paperabstract
}
\end{minipage}\vfill
}

\newcommand{\MakeFirstPageOneAuthor}{
\begin{center}
  \makepapertitle
  \paperfirstauthor
  \makeabstract
  \begin{minipage}[t]{0.3\textwidth}
   \raggedleft {\bf Date (final version):} 
  \end{minipage}\hspace{0.01\textwidth}
  \begin{minipage}[t]{0.6\textwidth}
    \today\\ 
    (submitted on March 31, 2021;\\ accepted on November 1, 2021)
  \end{minipage}\\[4mm]
\noindent%
  \begin{minipage}[t]{0.3\textwidth}
    \raggedleft {\bf Running head:} 
  \end{minipage}\hspace{0.01\textwidth}
  \begin{minipage}[t]{0.6\textwidth}
    \runninghead
  \end{minipage}\\[4mm]
\noindent%
  \begin{minipage}[t]{0.3\textwidth}
    \raggedleft {\bf How to cite:} 
  \end{minipage}\hspace{0.01\textwidth}
  \begin{minipage}[t]{0.6\textwidth}
    Linear Algebra Appl. {\bf 634} (2022), 179--209.\\[1mm]
   \footnotesize{\url{https://doi.org/10.1016/j.laa.2021.11.001}}
  \end{minipage}\\
\bigskip
\vfill
\end{center}

\renewcommand{\thefootnote}{}
\footnotetext[2]{2020 {\it Mathematics Subject Classification:} \thesubjclass}
\footnotetext[3]{{\it Key words and phrases:} \thekeywords}
\setcounter{footnote}{0}
\renewcommand{\thefootnote}{{\bf\,\alph{footnote}\alph{footnote}\alph{footnote}\,}}

\setcounter{page}{0}
\newpage
}

\makeatletter
\newcommand{\rightorleftmark}{%
  \begingroup\protected@edef\x{\rightmark}%
  \ifx\x\@empty
    \endgroup\leftmark
  \else
    \endgroup\rightmark
  \fi}
\makeatother

\newcommand{\papertitle}%
{Resolvent and spectrum for discrete symplectic systems in the limit point case}

\newcommand{\runninghead}%
{Resolvent, spectrum, and discrete symplectic systems}

\newcommand{\firstauthor}%
{Petr Zem{\'{a}}nek}
\newcommand{\firstauthorhead}%
{P. Zem{\'{a}}nek}

\newcommand{\firstaddress}
{Department of Mathematics and Statistics, Faculty of Science, Masaryk University \\
Kotl{\'{a}}{\v{r}}sk{\'{a}} 2, CZ-61137 Brno, Czech Republic}

\newcommand{\firstemail}%
{E-mail: zemanekp@math.muni.cz}

\newcommand{\paperabstract}%
{The spectrum of an arbitrary self-adjoint extension of the minimal linear relation associated with the discrete 
symplectic system in the limit point case is completely characterized by using the limiting Weyl--Titchmarsh 
\texorpdfstring{$M_+(\la)$}{M₊(λ)}-function. Furthermore, a dependence of the spectrum on a boundary condition is 
investigated and, consequently, several results of the singular Sturmian theory are derived.}

\newcommand{\thekeywords}%
{Discrete symplectic system; spectrum; eigenvalue; limit point case; \texorpdfstring{$M(\la)$}{M(λ)}-function.}

\newcommand{\thesubjclass}%
{{\it Primary\/} 47B39; {\it Secondary\/} 47A10; 47A06; 39A06; 39A12.}

\newcommand{\submittedto}%
{Linear Algebra and its Applications}

\DeclareMathAccent{\wwtilde}{\mathord}{largesymbols}{"65}
\DeclareMathSymbol{\widetildesym}{\mathord}{largesymbols}{"65}
\newcommand\lowerwidetildesym{%
  \text{\smash{\raisebox{-1.3ex}{%
    $\widetildesym$}}}}
\newcommand\wtilde[1]{%
  \mathchoice
    {\accentset{\displaystyle\lowerwidetildesym}{#1}}
    {\accentset{\textstyle\lowerwidetildesym}{#1}}
    {\accentset{\scriptstyle\lowerwidetildesym}{#1}}
    {\accentset{\scriptscriptstyle\lowerwidetildesym}{#1}}
}

\newcommand\lowerwidetildesymW{%
  \text{\smash{\raisebox{-1.35ex}{$\widetildesym$}}}}
\newcommand\wtildeW[1]{%
    {\accentset{\scalebox{1.2}{\lowerwidetildesymW}}{#1}}
}

\DeclareMathSymbol{\widehatsym}{\mathord}{largesymbols}{"62}
\newcommand\lowerwidehatsym{%
  \text{\smash{\raisebox{-1.3ex}{%
    $\widehatsym$}}}}
\newcommand\what[1]{%
  \mathchoice
    {\accentset{\displaystyle\lowerwidehatsym}{#1}}
    {\accentset{\textstyle\lowerwidehatsym}{#1}}
    {\accentset{\scriptstyle\lowerwidehatsym}{#1}}
    {\accentset{\scriptscriptstyle\lowerwidehatsym}{#1}}
}

\DeclareMathAlphabet{\mthdtcl}{U}{dutchcal}{m}{n}
\DeclareMathAlphabet{\mathpzc}{OT1}{pzc}{m}{it}
\DeclareMathAlphabet{\msfsl}{U}{eus}{m}{n}

\renewcommand{\d}{\mathrm{d}}
\newcommand{\dtau}{\d\tau}

\newcommand{\Jc}{\mathcal{J}}

\newcommand{\Oc}{\mathcal{O}}
\renewcommand{\O}{\mathrm{O}}

\newcommand{\Rc}{\mathcal{R}}
\newcommand{\Rm}{\mathrm{R}}
\newcommand{\Sc}{\mathcal{S}}

\newcommand{\Vc}{\mathcal{V}}

\newcommand{\Wc}{\mathcal{W}}
\newcommand{\tWc}{\wtildeW{\Wc}}
\newcommand{\Xc}{\mathcal{X}}

\newcommand{\mL}{\mathscr{L}}

\newcommand{\Cbb}{\mathbb{C}}

\newcommand{\Nbb}{\mathbb{N}}
\newcommand{\Rbb}{\mathbb{R}}
\newcommand{\Sbb}{\mathbb{S}}

\newcommand{\Zbb}{\mathbb{Z}}

\newcommand{\al}{\alpha}
\newcommand{\hal}{\what{\al}}
\newcommand{\be}{\beta}

\newcommand{\la}{\lambda}

\newcommand{\Ups}{\Upsilon}

\newcommand{\bla}{\bar{\la}}

\newcommand{\de}{\delta}
\newcommand{\De}{\Delta}

\newcommand{\Ps}{\Psi}

\newcommand{\Ph}{\Phi}

\newcommand{\eps}{\varepsilon}
\newcommand{\Ga}{\Gamma}

\newcommand{\ga}{\gamma}

\newcommand{\rh}{\rho}
\newcommand{\om}{\omega}
\newcommand{\tom}{\wtilde{\om}}
\newcommand{\si}{\sigma}

\newcommand{\sic}{\si_{\mathrm{c}}}
\newcommand{\sipc}{\si_{\mathrm{pc}}}
\newcommand{\sip}{\si_{\mathrm{d}}}
\newcommand{\sie}{\si_{\mathrm{e}}}
\newcommand{\TLP}{T_{\mathrm{LP}}}

\newcommand{\stm}{\hspace*{0.2mm}\fgebackslash\hspace*{0.3mm}}

\usepackage{scalerel}
\usepackage{stackengine}

\stackMath

\newcommand{\fclass}{[f]}

\newcommand{\gclass}{[g]}

\newcommand{\zclass}{[z]}
\newcommand{\tZ}{\wtilde{Z}}
\newcommand{\hz}{\hat{z}}
\newcommand{\hZ}{\what{Z}}

\newcommand{\ltp}{\ell^{\hspace{0.2mm}2}_{\Ps}}

\newcommand{\tltp}{\tilde{\ell}^{\hspace{0.3mm}2}_{\Ps}}

\newcommand{\tltpt}{\tilde{\ell}^{\hspace{0.3mm}2\times2}_{\Ps}}
\newcommand{\Tmax}{T_{\mathrm{max}}}
\newcommand{\Tmin}{T_{\mathrm{min}}}

\newcommand{\sZbb}{{\scriptscriptstyle{\Zbb}}}
\newcommand{\oinftyZ}{[0,\infty)_\sZbb}

\newcommand{\onZ}{[0,N]_\sZbb}

\newcommand{\mmatrix}[1]{\begin{pmatrix} #1
  \end{pmatrix}}
\newcommand{\msmatrix}[1]{\left(\begin{smallmatrix} #1
  \end{smallmatrix}\right)}  

\newcommand{\qtextq}[1]{\quad\text{#1}\quad}
\newcommand{\qtext}[1]{\quad\text{#1 }\ }

\DeclareMathOperator{\re}{Re}
\DeclareMathOperator{\im}{Im}

\DeclareMathOperator{\ran}{ran}
\DeclareMathOperator{\rank}{rank}
\DeclareMathOperator{\dom}{dom}

\DeclareMathOperator{\cotan}{cotan}

\usepackage{mathtools,xparse}
\renewcommand{\.}{\hspace*{0.1 em}}
\DeclarePairedDelimiter\xnorm{\lVert}{\rVert}
\NewDocumentCommand{\norm}{som}
 {\IfBooleanTF{#1}
   {\xnorm*{#3}}
   {\IfNoValueTF{#2}
     {\mathopen{|\mkern-.8mu|}\.#3\.\mathclose{|\mkern-.8mu|}}
     {\xnorm[#2]{\.#3\.}}%
   }
 }
\DeclarePairedDelimiter\xinner{\langle}{\rangle}

\NewDocumentCommand{\xinnr}{som}
 {\IfBooleanTF{#1}
   {\xinner*{#3}}
   {\IfNoValueTF{#2}
     {\mathopen{\langle}\.#3\.\mathclose{\rangle}}
     {\xinner[#2]{\.#3\.}}%
   }
 }

\makeatletter
\def\inner{\@ifnextchar[{\@INNwith}{\@INNwithout}}
\def\@INNwith[#1]#2#3{\xinnr[#1]{#2,#3}}
\def\@INNwithout#1#2{\xinnr{#1,#2}}

\def\innerP{\@ifnextchar[{\@INNPwith}{\@INNPwithout}}
\def\@INNPwith[#1]#2#3{\xinnr[#1]{#2,#3}_\Ps}
\def\@INNPwithout#1#2{\xinnr{#1,#2}_\Ps}

\def\innerPN{\@ifnextchar[{\@INNPNwith}{\@INNPNwithout}}
\def\@INNPNwith[#1]#2#3{\xinnr[#1]{#2,#3}_{\Ps,N}}
\def\@INNPNwithout#1#2{\xinnr{#1,#2}_{\Ps,N}}

\def\abs{\@ifnextchar[{\@awith}{\@awithout}}
\def\@awith[#1]#2{{#1|}#2\.{#1|}}
\def\@awithout#1{|#1\.|}
 
\def\normS{\@ifnextchar[{\@Nwith}{\@Nwithout}}
\def\@Nwith[#1]#2{\norm[#1]{#2}_\si}
\def\@Nwithout#1{\norm{#1}_\si}

\def\normE{\@ifnextchar[{\@NEwith}{\@NEwithout}}
\def\@NEwith[#1]#2{\norm[#1]{#2}_2}
\def\@NEwithout#1{\norm{#1}_2}

\def\normA{\@ifnextchar[{\@NAwith}{\@NAwithout}}
\def\@NAwith[#1]#2{\norm[#1]{#2}_1}
\def\@NAwithout#1{\norm{#1}_1}

\def\normP{\@ifnextchar[{\@NPwith}{\@NPwithout}}
\def\@NPwith[#1]#2{\norm[#1]{#2}_{\Ps}}
\def\@NPwithout#1{\norm{#1}_{\Ps}}

\def\normW{\@ifnextchar[{\@NWwith}{\@NWwithout}}
\def\@NWwith[#1]#2{\norm[#1]{#2}_{\Wc}}
\def\@NWwithout#1{\norm{#1}_{\Wc}}

\def\normtW{\@ifnextchar[{\@NtWwith}{\@NtWwithout}}
\def\@NtWwith[#1]#2{\norm[#1]{#2}_{\scalebox{0.7}{$\tWc$}}}
\def\@NtWwithout#1{\norm{#1}_{\scalebox{0.6}{$\tWc$}}}

\newcommand{\Sla}[1]{\text{\rm(S$_{#1}$})}
\newcommand{\Slaf}[2]{\text{\rm(S$_{#1}^{#2}$)}}
\DeclareMathOperator{\diag}{diag}
\DeclareMathOperator{\sgn}{sgn}

\makeatletter
\newcommand{\ltxlabel}{\ltx@label}
\makeatother
\newcounter{GatherItemCounter}

\usepackage{ifdraft}
\ifdraft{\overfullrule26pt}{\overfullrule0pt}
\ifdraft{\usepackage[nomsgs]{refcheck}}{}
\ifdraft{\geometry{showframe,showcrop}}{}

\hypersetup{
pdfencoding=unicode,
pdftitle={\papertitle},
pdfauthor={\firstauthor},
pdfsubject={Department of Mathematics and Statistics, Faculty of Science, Masaryk University},
pdfkeywords={\thekeywords},
pdfsubject=\paperabstract,
pdfencoding=auto
}

\begin{document}

%%%%%%%%%%%%%%%%%%%%%%%%%%%%%%%%%%%%%%%%%%%%%% FIRST PAGE %%%%%%%%%%%%%%%%%%%%%%%%%%%%%%%%%%%%%%%%%%%%%%%%%%%%%%%%%%%%

\MakeFirstPageOneAuthor

% \tableofcontents

%%%%%%%%%%%%%%%%%%%%%%%%%%%%%%%%%%%%%%%%%%%% SECTION %%%%%%%%%%%%%%%%%%%%%%%%%%%%%%%%%%%%%%%%%%%%%%%%%%%%%%%%%%%%%%%%%

\section{Introduction}\label{S:intro}

{\it Every operator has a secret treasure: its spectrum} (A.~Khrabustovskyi, O.~Post, and C.~Trunk in the summary of 
their special session at IWOTA 2019). Motivated by this quote, we investigate the spectrum of ``operators'', which are 
beyond a certain class of the {\it time-reversed discrete symplectic systems} on the half-line, i.e.,
 \begin{equation*}\label{E:Sla}\tag{S$_\la$}
  z_{k}(\la)=(\Sc_k+\la\.\Vc_k)\,z_{k+1}(\la), \quad k\in\oinftyZ\coloneq [0,\infty)\cap\Zbb,
 \end{equation*}
with $\la\in\Cbb$ playing the role of the spectral parameter and the coefficients being $2n\times2n$ 
complex-valued matrices such that 
 \begin{equation*}
  \Sc_k^*\Jc\Sc_k=\Jc, \quad 
  \Vc_k^*\Jc\Sc_k\ \text{ is Hermitian,} \qtextq{and} 
  \Vc_k^*\Jc\,\Vc_k=0 \qtext{for all $k\in\oinftyZ$,}
 \end{equation*}
where the superscript $*$ denotes the conjugate transpose and $\Jc$ stands for the $2n\times2n$ orthogonal and 
skew-symmetric matrix 
 \begin{equation*}
  \Jc=\Jc_{2n}\coloneq\mmatrix{0 & I_n\\ -I_n & 0}.
 \end{equation*}
Since system~\eqref{E:Sla} can be equivalently written as $\Jc(z_k-\Sc_k\,z_{k+1})=\la\.\Ps_k\.z_k$ with 
$\Ps_k\coloneq \Jc\.\Sc_k\.\Jc\.\Vc^*_k\.\Jc$, it gives rise to the linear mapping from the space of $2n$-vector valued 
square summable sequences on $\oinftyZ$ to itself given by
 \begin{equation*}%\label{E:mL.def}
  \mL(z)_k\coloneq\Jc(z_k-\Sc_k\,z_{k+1}),
 \end{equation*}
for which the square summability is defined via the semi-inner product 
 \begin{equation*}%\label{E:semi-inner.def}
  \innerP{z}{u}\coloneq\sum_{k=0}^\infty z_k^*\,\Ps_k\,u_k \qtextq{and the induced semi-norm}
  \normP{z}\coloneq\sqrt{\innerP{z}{z}}
 \end{equation*} 
with respect to the weight matrices $\Ps_k\geq0$, i.e., we restrict our attention to the space 
 \begin{equation*}%\label{E:ltp.def}
  \ltp=\ltp(\oinftyZ)\coloneq\{z\in\Cbb(\oinftyZ)^{2n}\mid \normP{z}<\infty\}.
 \end{equation*} 
However, according to \cite[Corollary~1]{pZ20}, this natural mapping does yield only a non-densely defined operator on 
the Hilbert space 
 \begin{equation*}%\label{E:tltp.def}
  \tltp=\tltp(\oinftyZ)\coloneq \ltp\big/\big\{z\in\Cbb(\oinftyZ)^{2n}\mid \ \normP{z}=0\big\}
 \end{equation*} 
consisting of equivalence classes denoted by the square brackets as $\zclass$ and with the inner product 
$\innerP{\zclass}{\fclass}\coloneq\innerP{z}{f}$ for any $\zclass,\fclass\in\tltp$. This fact forces the approach based 
on linear relations instead of operators as it was already done by the author and his collaborator 
in~\cite{slC.pZ15,pZ.slC16,pZ.slC:SAE2}, where the main attention was paid to various characterizations of self-adjoint 
extensions of the minimal linear relation $\Tmin\subseteq\Tmax\subseteq\tltpt\coloneq\tltp\times\tltp$ being the 
adjoint of the maximal linear relation, i.e., it is the closed and symmetric relation such that $\Tmin=\Tmax^*$ for
 \begin{align*}
  &\Tmax\coloneq \big\{\{\zclass,\fclass\}\in\tltpt\mid \text{there exists $u\in\zclass$ such that }
                                                        \mL(u)_k=\Ps_k\.f_k\ \text{ for all } k\in\oinftyZ\big\}.
 \end{align*}
In the present manuscript, we focus on the spectrum of an arbitrary self-adjoint extension $\TLP$ of $\Tmin$ related to 
system~\eqref{E:Sla} being in the limit point case (i.e., it possesses precisely $n$ linearly independent square 
summable solutions for all $\la\in\Cbb\stm\Rbb$) and satisfying the {\it strong Atkinson condition} (i.e., 
$\normP{z}>0$ for every nontrivial solution of the system), which guarantees that for any pair 
$\{\zclass,\fclass\}\in\Tmax$ there is a unique representative $\hz\in\zclass$ with the property $\mL(\hz)=\Ps\.f$ on 
$\oinftyZ$. By \cite[Corollary~3.5]{pZ.slC16} or \cite[Corollary~4.7]{pZ.slC:SAE2}, this linear 
relation admits the representation
 \begin{equation*}
  \TLP=\TLP(\al)=\big\{\{\zclass,\fclass\}\in\Tmax\mid \al\.\hz_0=0\big\}
 \end{equation*}
for a suitable 
 \begin{equation*}
  \al\in\Ga\coloneq\{\al\in\Cbb^{n\times2n}\mid \al\.\al^*=I,\ \al\.\Jc\al^*=0\}.
 \end{equation*}
Since the relevant part of the abstract theory of linear relations was recently summarized in~\cite{pZ.slC:SAE2} and a 
much more comprehensive treatise can be found, e.g., in~\cite{jB.sH.hsvdS20}, we recall now only that the {\it 
spectrum} $\si(\TLP)$ is the complement of the {\it resolvent set} $\rh(\TLP)$ consisting of all $\la\in\Cbb$ such that 
$(\TLP-\la\.I)^{-1}$ is the graph of a certain everywhere defined bounded operator on $\tltp$, in which case 
$\Rc_{\TLP}(\la)\coloneq(\TLP-\la\.I)^{-1}$ represents the {\it resolvent relation}. Then $\si(\TLP)\subseteq\Rbb$, it 
is a closed set and it can be expressed as the union of three disjoint sets 
$\si(\TLP)=\sip(\TLP)\cup\sipc(\TLP)\cup\sic(\TLP)$, where
 \begin{enumerate}[leftmargin=10mm,topsep=2mm,label={{\normalfont{(\roman*)}}}]
  \item $\sip(\TLP)$ denotes the {\it discrete spectrum} consisting of isolated points of the spectrum, while all other 
        points (i.e., accumulation points) of the spectrum belong to the {\it essential spectrum} 
        $\sie(\TLP)\coloneq \si(\TLP)\stm\sip(\TLP)=\sipc(\TLP)\cup\sic(\TLP)$;
  \item eigenvalues in $\sie(\TLP)$ form the {\it point-continuous spectrum} $\sipc(\TLP)$;
  \item and $\sic(\TLP)\coloneq \sie(\TLP)\stm\sipc(\TLP)$ stands for the {\it continuous spectrum}.
 \end{enumerate}
Let us emphasize that $\sip(\TLP)$ coincides with the set of all isolated eigenvalues determined by 
system~\eqref{E:Sla} and the boundary condition $\al\.z_0(\la)=0$ or equivalently such that 
$\ker(\TLP-\la\.I)\neq\{0\}$. In particular, if $\si(\TLP)=\sip(\TLP)$, then $\TLP$ is said to have {\it pure discrete 
spectrum}, while for $\si(\TLP)=\sic(\TLP)$ we have {\it pure continuous spectrum}. The essential spectrum could also 
include eigenvalues of infinite multiplicity, but it is not possible in the present case as the multiplicity of $\la$ 
is 
equal to $\dim\ker(\TLP-\la\.I)$. For completeness, let us note that another approach to the structure of the spectrum 
can be found in~\cite[Definition~1.2.3]{jB.sH.hsvdS20}, but for the linear relation $\TLP$ the above described 
structure 
of the spectrum is fully sufficient and it is the same as for operators considered, e.g., 
in~\cite{jC.wnE68,dbH.jkS82:QM,yS06}. 

The limit point case hypothesis is crucial for the present research, because it guarantees the uniqueness of the 
definition of the $n\times n$ matrix-valued {\it limiting Weyl--Titchmarsh function} as
 \begin{equation*}%\label{E:Mla.def}
  M_+(\la)\coloneq \lim_{N\to\infty} M_N(\la,\al,\be)
  \qtextq{with}
  M_N(\la,\al,\be)\coloneq-[\be\.\tZ_N(\la)]^{-1}\.\be\.\hZ_N(\la) 
 \end{equation*}
through any $\be\in\Ga$ and the pair of $2n\times n$ matrix-valued solutions $\hZ(\la)$ and $\tZ(\la)$ determined by 
the initial conditions 
 \begin{equation}\label{E:hz.tz.def}
  \hZ_0(\la)=\al^* \qtextq{and} \tZ_0(\la)=-\Jc\al^*.
 \end{equation}
Indeed, for all $N$ large enough, every $M_{N+1}(\la,\al,\be)$ belongs to the boundary of a~certain Weyl disks; these 
disks are nested as $N$ grows and they collapse in the limit point case to the same singleton for all $\be\in\Ga$, so 
the limit in the definition of $M_+(\la)$ above truly does not depend on the choice of $\be$ and it coincides with the 
center of the so-called limiting Weyl disk. Consequently, the columns of the corresponding $2n\times n$ matrix-valued 
{\it Weyl solution}
 \begin{equation*}
  \Xc^+(\la)\coloneq \hZ(\la)+\tZ(\la)\.M_+(\la)
 \end{equation*}
span all linearly independent square summable solutions of~\eqref{E:Sla} for all $\la\in\Cbb\stm\Rbb$, see 
\cite{slC.pZ10,rSH.pZ14:JDEA} for more details. The $M_+(\la)$-function is well defined for all 
$\la\in\Cbb\stm\Rbb$ and its behavior for $\la\in\Rbb$ is closely related to the spectrum of $\TLP$ as we can see from 
our main result stated below. 

\begin{theorem}\label{T:intro.spectrum}
 Let $\al\in\Ga$ be given, system~\eqref{E:Sla} be in the limit point case for all $\la\in\Cbb\stm\Rbb$, and the strong 
 Atkinson condition hold.
 \begin{enumerate}[leftmargin=10mm,topsep=1mm,label={{\normalfont{(\roman*)}}}]
  \item We have $\la_0\in\rh(\TLP)$ if and only if $M_+(\la)$ is holomorphic at $\la_0$. Furthermore, in this case the 
        resolvent relation admits the representation
         \begin{equation*}%\label{E:MR1}
          \Rc_{\TLP}(\la_0)\coloneq(\TLP-\la_0\.I)^{-1}
           =\Bigg\{\bigg\{\fclass,\Big[\sum_{j=0}^\infty G_{kj}(\la_0)\.\Ps_j\.f_j\Big]\bigg\}\,\Big|\, f\in\ltp\Bigg\},
         \end{equation*}
        where
         \begin{equation*}%\label{E:MR2}
          G_{kj}(\la_0)\coloneq
           \begin{cases}
            \tZ_k(\la_0)\.\Xc^{+*}_j(\bla_0),& k\in[0,j]_\sZbb,\\
            \Xc^+_k(\la_0)\.\tZ^*_j(\bla_0),& k\in[j+1,\infty)_\sZbb.
           \end{cases}
         \end{equation*}
  
  \item A number $\la_0\in\Rbb$ belongs to $\sip(\TLP)$, i.e., it is an isolated eigenvalue of $\TLP$, if and only if 
        $M_+(\la)$ has a simple pole at $\la_0$, which is equivalent to the existence of the Laurent expansion of 
        $M_+(\la)$ in a neighborhood of $\la_0$ as
         \begin{equation*}%\label{E:MR5}
          M_+(\la)=K_{-1}\.(\la-\la_0)^{-1}+K_0+K_1\.(\la-\la_0)+\cdots,
         \end{equation*}
        where $K_{-1},K_0,\dots$ are $n\times n$ Hermitian matrices with
         \begin{equation*}%\label{E:MR6}
          K_{-1}=-\big(\tau(\la_0)-\tau(\la_0^-)\big)\lneqq0,
         \end{equation*}
        where $\tau$ stands for the corresponding limiting spectral function. In this case, all columns of 
        $\tZ(\la_0)\.K_{-1}$ belong to $\ltp$ and the nonzero columns are eigenfunctions of $\TLP$ related to $\la_0$. 
        Furthermore, the columns of 
         \begin{equation*}
          \hZ(\la_0)+\tZ(\la_0)\.K_0+\Big[\frac{\d}{\d\la} \tZ(\la)\Big]_{\la=\la_0}\.K_{-1}
         \end{equation*}
        belong also to $\ltp$. 
  
  \item We have $\la_0\in\sipc(\TLP)$ if and only if $M_+(\la)$ is not holomorphic at $\la_0$, it holds
         \begin{equation*}%\label{E:MR26}
          L\coloneq \lim_{\nu\to0} \nu\.M_+(\la_0+i\.\nu)\neq0,
         \end{equation*}
        and $M_+(\la)-i\.L\.(\la-\la_0)^{-1}$ is not holomorphic at $\la_0$. In this case, the columns of 
        $\tZ(\la_0)\.L$ belong to $\ltp$ and its nonzero columns are eigenfunctions of $\TLP$ related to $\la_0$.

  \item We have $\la_0\in\sic(\TLP)$ if and only if $M_+(\la)$ is not holomorphic at $\la_0$ and 
         \begin{equation*}
          \lim_{\nu\to0} \nu\.M_+(\la_0+i\.\nu)=0.
         \end{equation*}
 \end{enumerate}
\end{theorem}

The proof of the latter theorem is given in Section~\ref{S:spectrum} in four separate statements, see 
Theorems~\ref{T:MR1}--\ref{T:MR4}. Our approach is based on a generalization of some ideas from the continuous case 
studied in~\cite{dbH.jkS82:QM}, e.g., the expansion of $M_+(\la)$-function or a connection between its poles and jumps 
of the limiting spectral function $\tau$ derived in Section~\ref{S:prelim}. However, their use in the context of 
discrete symplectic systems forces one to work with linear relations. Furthermore, our results may provide an 
alternative way for the development of the so-called {\it singular Sturmian theory} concerning, e.g., the dependence of 
eigenvalues on the choice of the matrix $\al$, which seems to be not completely solved for discrete symplectic systems, 
see Subsection~4.2.4 of the encyclopedia of discrete symplectic systems \cite{oD.jE.rSH19}. These applications are 
discussed in Section~\ref{S:application}.

Finally, it should be mentioned that discrete symplectic systems naturally arise in the discrete calculus of 
variations, numerical integration schemes or in the theory of continued fractions. They include any even-order 
Sturm--Liouville difference equation and linear Hamiltonian difference system and system~\eqref{E:Sla} is the proper 
discrete analogue of the {\it linear Hamiltonian} (or {\it canonical}) {\it differential systems}
 \begin{equation}\label{E:cLHS}\tag{H$_\la$}
  -\Jc\.z'(t,\la)=[H(t)+\la\.W(t)]\.z(t,\la),\quad t\in(a,b),
 \end{equation}
where $H(t)$ and $W(t)\geq0$ are $2n\times 2n$ Hermitian-matrix valued functions, see e.g. 
\cite{slC.pZ10,rSH.pZ14:AMC,rSH.pZ14:ICDEA}. The latter fact justifies our effort in the development of the spectral 
theory for these systems. The present results can be seen as a discrete counterpart of~\cite{dbH.jkS82:QM} for 
system~\eqref{E:cLHS} and as a generalization of a~correction of~\cite[Sections~6 and~7]{yS06} for a linear Hamiltonian 
difference system, where the non-density of the domain of the related minimal operator is omitted; see also 
\cite{foN17,jC.wnE68,hS.yS14,rB.jE.aK.gT14,hS.qK.yS16,joA.doA.foN13,hB13,hS.yS15}. 

%%%%%%%%%%%%%%%%%%%%%%%%%%%%%%%%%%%%%%%%%%%% SECTION %%%%%%%%%%%%%%%%%%%%%%%%%%%%%%%%%%%%%%%%%%%%%%%%%%%%%%%%%%%%%%%%%
 
\section{Preliminaries}\label{S:prelim}

For reader's convenience, let us recall the notation used throughout the paper. The real and imaginary parts 
of any $\la\in\Cbb$ are, respectively, denoted by $\re\la$ and $\im\la$. The symbols $\Cbb_+$ and $\Cbb_-$ mean, 
respectively, the upper and lower complex planes, i.e., $\Cbb_+\coloneq\{\la\in\Cbb\mid \im(\la)>0\}$ and 
$\Cbb_-\coloneq\{\la\in\Cbb\mid \im(\la)<0\}$. All matrices are considered over the field of complex numbers $\Cbb$. 
For 
$r,s\in\Nbb$ we denote by $\Cbb^{r\times s}$ the space of all complex-valued $r\times s$ matrices and $\Cbb^{r\times 
1}$ 
will be abbreviated as $\Cbb^r$. In particular, the $r\times r$ \emph{identity} and \emph{zero matrices} are written as 
$I_r$ and $0_r$, where the subscript is omitted whenever it is not misleading (for simplicity, the zero vector is also 
written as $0$). For a given matrix $M\in\Cbb^{r\times s}$ or $M\in\Cbb^{r\times r}$ we indicate by $M^*$, $\ker M$, 
$\ran M$, $\det M$, $\rank M$, $\re M\coloneq (M+M^*)/2$, $\im M\coloneq (M-M^*)/(2i)$, $M\geq0$, and $M>0$ 
respectively, its conjugate transpose, kernel, range (or image), determinant, rank, Hermitian components (real and 
imaginary parts), positive semidefiniteness, and positive definiteness. Finally, the symbol $\diag\{M_1,\dots,M_n\}$ 
stands for the block-diagonal matrix with matrices $M_1,\dots,M_n$ located on the main diagonal.

By $\Cbb(\oinftyZ)^{r\times s}$ we denote the space of sequences defined on $\oinftyZ$ of complex $r\times s$ matrices, 
where typically $r\in\{n,2n\}$ and $1\leq s\leq2n$. In particular, we write only $\Cbb(\oinftyZ)^r$ in the case $s=1$. 
If $M\in\Cbb(\oinftyZ)^{r\times s}$, then $M(k)\coloneq M_k$ for $k\in\oinftyZ$ and if 
$M(\la)\in\Cbb(\oinftyZ)^{r\times s}$, then $M(\la,k)\coloneq M_k(\la)$ for $k\in\oinftyZ$ with 
$M_k^*(\la)\coloneq[M_k(\la)]^*$. If $M\in\Cbb(\oinftyZ)^{r\times s}$ and $N\in\Cbb(\oinftyZ)^{s\times p}$, then 
$MN\in\Cbb(\oinftyZ)^{r\times p}$, where $(MN)_k\coloneq M_k N_k$ for $k\in\oinftyZ$. Finally, the forward difference 
operator acting on $\Cbb(\oinftyZ)^{r\times s}$ is denoted by $\De$ with $(\De z)_k\coloneq\De z_k$ 
and we put $z_k \big|_{m}^n\coloneq z_n - z_m$.

In the following hypothesis we summarize the basic assumptions concerning system~\eqref{E:Sla} or its 
nonhomogeneous counterpart
 \begin{equation}\label{E:Slaf}\tag{S$_\la^f$}
  z_{k}(\la)=(\Sc_k+\la\.\Vc_k)\,z_{k+1}(\la)-\Jc\,\Ps_k\,f_k, \quad k\in\oinftyZ.
 \end{equation}
Let us emphasize that these systems can be determined either by the pair of coefficient matrices 
$\{\Sc,\Vc\}$ or by the pair $\{\Sc,\Ps\}$ and there is no difference between these two approaches as the matrices 
$\Ps_k$ and $\Vc_k$ are mutually connected via the equalities $\Ps_k=\Jc\.\Sc_k\.\Jc\.\Vc^*_k\.\Jc$ and 
$\Vc_k=-\Jc\,\Ps_k\,\Sc_k$. 

\begin{remark}
 The last part of the hypothesis concerns the equivalent formulation of the strong Atkinson 
 condition, which is a traditional assumption in the Weyl--Titchmarsh theory for differential or difference equations 
 and we utilize it because of the uniqueness property mentioned in the introduction. However, it is worth noticing 
 that it is not needed to require the property $\normP{z(\la)}>0$ for all $\la\in\Cbb$ explicitly, see 
 \cite[Lemma~2.2]{pZ?:MN}. Furthermore, let us point out that Theorem~\ref{T:Lagrange} does not rely on the limit point 
 case and for Theorem~\ref{T:A3} the weaker version of this assumption would be sufficient, see 
 also~\cite{rSH.pZ14:JDEA}.
\end{remark}
 
\begin{hypothesis}\label{H:basic}
 A number $n\in\Nbb$ and matrix $\al\in\Ga$ are given and we have the pair of matrix-valued sequences 
 $\Sc,\Ps\in\Cbb(\oinftyZ)^{2n\times2n}$ such that
  \begin{equation*}%\label{E:H.basic}
   \Sc_k^*\Jc\Sc_k=\Jc,\quad \Ps_k^*=\Ps_k, \quad \Ps_k^*\,\Jc\,\Ps_k=0, \qtextq{and} \Ps_k\geq0 
   \qtext{for all $k\in\oinftyZ$.}
  \end{equation*}
 We also define $\Vc_k\coloneq -\Jc\,\Ps_k\,\Sc_k$ for all $k\in\oinftyZ$. Furthermore, system~\eqref{E:Sla} is in the 
 limit point case for all $\la\in\Cbb\stm\Rbb$ (i.e., it possesses $n$ linearly independent square summable solutions) 
 and we assume that there are $\la\in\Cbb$ and $N_0\in\oinftyZ$ such that every nontrivial solution of~\eqref{E:Sla} 
 satisfies $\sum_{k=0}^{N_0} z_k^*(\la)\.\Ps_k\.z_k(\la)>0$.
\end{hypothesis}

The following identity represents one of the most important tools in the whole theory related to system~\eqref{E:Sla}. 
For its proof we refer to~\cite[Theorem~2.5]{slC.pZ15}.

\begin{theorem}[Extended Lagrange formula]\label{T:Lagrange}
 Let Hypothesis~\ref{H:basic} be satisfied, numbers $\la,\nu\in\Cbb$, $m\in\{1,\dots,2n\}$ and sequences 
 $f,g\in\Cbb(\oinftyZ)^{2n\times m}$ be given. If $z(\la),u(\nu)\in\Cbb(\oinftyZ)^{2n\times m}$ are solutions of 
 systems~\eqref{E:Slaf} and \Slaf{\nu}{g}, respectively, then for any $k,s,t\in\oinftyZ$ such that $s\leq t$, we have 
  \begin{align*}
   \De[z_k^*(\la)\,\Jc u_k(\nu)]&=(\bla-\nu)\,z_k^*(\la)\,\Ps_k\,u_k(\nu)+f_k^*\,\Ps_k\,u_k(\nu)
                                                            -z_k^*(\la)\,\Ps_k\,g_k,\notag\\%\label{E:2.6}\\
   &\hspace*{-33mm} z_k^*(\la)\,\Jc u_k(\nu)\big|_{s}^{t+1}
     =\sum_{k=s}^{t}\big\{(\bla-\nu)\,z_k^*(\la)\,\Ps_k\,u_k(\nu)+f_k^*\,\Ps_k\,u_k(\nu)
                                                            -z_k^*(\la)\,\Ps_k\,g_k\big\}.%\label{E:lagrange}
  \end{align*}
 Especially, if $\nu=\bla$ and $f\equiv0\equiv g$, then we get the Wronskian-type identity
  \begin{equation*}%\label{E:wronski.id}
   z_k^*(\la)\,\Jc\,u_k(\bla)\equiv z_0^*(\la)\,\Jc\,u_0(\bla) \qtext{on $\oinftyZ$.}
  \end{equation*}
\end{theorem} 

The existence of a unique solution $z(\la)\in\Cbb(\oinftyZ)^{2n}$ of any initial value problem associated with 
system~\eqref{E:Sla} is easily guaranteed by the symplectic-type equality 
 \begin{equation*}
  (\Sc_k+\la\.\Vc_k)^*\.\Jc\.(\Sc_k+\la\.\Vc_k)=\Jc
 \end{equation*}
for all $k\in\oinftyZ$ and $\la\in\Cbb$, which yields 
$\abs{\det(\Sc_k+\la\.\Vc_k)}\equiv1$. Furthermore, in~\cite[Theorem~4.3 and Corollary~4.4]{slC.pZ15}, we proved that 
for any $\la\in\Cbb\stm\Rbb$, $f\in\ltp$, and $\xi\in\Cbb^n$ the unique solution $\hz\in\ltp$ of the boundary value 
problem
 \begin{equation}\label{E:bvp}
  \eqref{E:Slaf} \qtextq{and} \al\.\hz_0=\xi
 \end{equation}
can be expressed as 
  \begin{equation*}
   \hz_k\coloneq \Xc^+_k(\la)\.\xi+\Rm_\la(f)_k,\quad k\in\oinftyZ,
  \end{equation*}
where $\Rm_\la(f)_k\coloneq\sum_{j=0}^\infty G_{kj}\.\Ps_j\.f_j$ can be seen a representative of the value of the 
resolvent relation $\Rc_{\TLP}(\la)$ and $G_{kj}(\la)$ stands for the {\it Green's function}
 \begin{equation*}
  G_{kj}(\la)\coloneq
   \begin{cases}
    \tZ_k(\la)\.\Xc^{+*}_j(\bla),& k\in[0,j]_\sZbb,\\
    \Xc^+_k(\la)\.\tZ^*_j(\bla),& k\in[j+1,\infty)_\sZbb.
   \end{cases}
 \end{equation*}

A necessary and sufficient condition for system~\eqref{E:Sla} being in the limit point case for all $\la\in\Cbb\stm\Rbb$
was established, e.g., in~\cite[Corollary~4.14]{rSH.pZ14:JDEA}. More specifically, this case occurs if and only if for 
every $\nu,\si\in\Cbb\stm\Rbb$ and every square summable solution $z(\nu)$ and $z(\si)$ of systems~\Sla{\nu} 
and~\Sla{\si}, respectively, we have
 \begin{equation}\label{E:lpc.equiv}
  \lim_{k\to\infty} z_k^*(\nu)\.\Jc\.z_k(\si)=0.
 \end{equation}
Equality~\eqref{E:lpc.equiv} is true also when $\nu,\si\in\Rbb$. But it does not follow from the latter statement, 
because its proof is based on a certain limit behavior of $M_+(\la)$-function, which may not exist for some real values 
of the parameter. Actually, the existence of $M_+(\la)$-function for real values of $\la$ is the main object of the 
present research. So, instead of a straightforward modification, we incorporate a non-real valued parameter in 
the system, which gives us the following extension. This result enables us to assume the limit point case only, which 
brings us back to~\cite[Equation~(3.3)]{jC.wnE68} in contrast to the approach utilized in~\cite{dbH.jkS82:QM}, where 
equality~\eqref{E:lpc.equiv} was required explicitly. For completeness, we note that another sufficient condition for 
the limit point case can be found in~\cite[Theorem~4.4]{pZ?:MN} and~\cite[Theorem~2.7]{pZ.slC16}.

\begin{theorem}\label{T:A3}
 Let Hypothesis~\ref{H:basic} hold. Then system~\eqref{E:Sla} is in the limit point case for all $\la\in\Cbb\stm\Rbb$ 
 if and only if for any sequences $z,y\in\ltp$ satisfying nonhomogeneous systems~\Slaf{0}{f} and~\Slaf{0}{h} for some 
 $f,h\in\ltp$, respectively, it holds
  \begin{equation}\label{E:P.A3.1}
   \lim_{k\to\infty} z_k^*\.\Jc\.y_k=0.
  \end{equation}
\end{theorem}
\begin{proof}
 Necessity is an immediate consequence of the previous discussion concerning 
 equality~\eqref{E:lpc.equiv}. More precisely, if~\eqref{E:P.A3.1} holds for any pair of square sum\-ma\-ble solutions 
 of~\Slaf{0}{f} and~\Slaf{0}{h}, then it is true also for any pair of square summable solutions of systems~\Sla{\nu} 
 and~\Sla{\si} with arbitrary $\nu,\si\in\Cbb\stm\Rbb$ as they correspond to the nonhomogeneous systems with 
 $f\coloneq \nu\.z$ and $h\coloneq\si\.y$, respectively. Consequently, system~\eqref{E:Sla} is in the limit point case 
 for all $\la\in\Cbb\stm\Rbb$.
 
 To show sufficiency, let system~\eqref{E:Sla} be in the limit point case for any $\la\in\Cbb\stm\Rbb$ and 
 $z,y\in\ltp$ be, respectively, solutions of~\Slaf{0}{f} and~\Slaf{0}{h} for some $f,h\in\ltp$. If we fix 
 $\la\in\Cbb\stm\Rbb$, then these two sequences solve also~\Slaf{\la}{F} and~\Slaf{\bla}{H}, respectively, with 
 \begin{equation*}
  F\coloneq f-\la\.z\in\ltp \qtextq{and} H\coloneq h-\bla\.y\in\ltp.
 \end{equation*}
 If we denote
  \begin{equation*}
   \xi_F\coloneq \al\.z_0 \qtextq{and} \xi_H\coloneq\al\.y_0,
  \end{equation*}
 then the sequences $z,y$ solve similar boundary value problems as in~\eqref{E:bvp}, which implies that
  \begin{equation}\label{E:P.A3.z}
   z_k=z_k(\la)=\Xc^+_k(\la)\.\xi_F+\Rm_\la(F)_k \qtext{for all $k\in\oinftyZ$}
  \end{equation}
 and
  \begin{equation}\label{E:P.A3.y}
   y_k=y_k(\bla)=\Xc^+_k(\bla)\.\xi_H+\Rm_{\bla}(H)_k \qtext{for all $k\in\oinftyZ$.}
  \end{equation}
 Upon applying Theorem~\ref{T:Lagrange} we obtain
  \begin{equation}\label{E:P.A3.3}
   \begin{aligned}
    \lim_{k\to\infty} z_k^*\.\Jc\.y_k-z_0^*\.\Jc\.y_0
     &=\innerP{F}{y}-\innerP{z}{H}\\
     &=\innerP{F}{\Xc^+(\bla)\.\xi_H}-\innerP{\Xc^+(\la)\.\xi_F}{H}\\
      &\hspace*{20mm}+\innerP{F}{\Rm_{\bla}(H)}-\innerP{\Rm_{\la}(F)}{H}.
   \end{aligned}
  \end{equation}
 By using the expressions of $z_0$ and $y_0$ from \eqref{E:P.A3.z}--\eqref{E:P.A3.y}, the special form of 
 the initial value of the fundamental matrix $\Ph(\la)\coloneq\big(\hZ(\la),\,\, \tZ(\la)\big)$ at $k=0$, 
 and the Hermitian property of $M_+(\la)$ we can verify that 
  \begin{equation*}
   z_0^*\.\Jc\.y_0=\innerP{\Xc^+(\la)\.\xi_F}{H}-\innerP{F}{\Xc^+(\bla)\.\xi_H}.
  \end{equation*}
 Simultaneously, the absolute convergence of the double sum 
  \begin{equation*}
   \innerP{F}{\Rm_{\bla}(H)}=\sum_{k=0}^\infty\sum_{j=0}^\infty F_k^*\.\Ps_k\.G_{kj}(\bla)\.\Ps_j\.H_j
  \end{equation*}
 guarantees that we can interchange the order of the summation, which yields
  \begin{align*}
   \innerP{F}{\Rm_{\bla}(H)}
    &=\sum_{j=0}^\infty\sum_{k=0}^\infty F_k^*\.\Ps_k\.G_{kj}(\bla)\.\Ps_j\.H_j
     =\sum_{j=0}^\infty\Big(\sum_{k=0}^{j-1}+\sum_{k=j}^{j}+\sum_{k=j+1}^\infty\Big) 
                                                            F_k^*\.\Ps_k\.G_{kj}(\bla)\.\Ps_j\.H_j\\
    &=\sum_{j=0}^\infty\bigg\{\Big(\sum_{k=0}^{j-1}+\sum_{k=j+1}^\infty\Big) F_k^*\.\Ps_k\.G^*_{jk}(\la)\.\Ps_j\.H_j 
      +F_j^*\.\Ps_j\.(G_{jj}^*(\la)-\Jc)\.\Ps_j\.H_j\bigg\}\\
    &=\sum_{j=0}^\infty\Big\{\sum_{k=0}^\infty F_k^*\.\Ps_k\.G^*_{jk}(\la)\.\Ps_j\.H_j
          -F_j^*\.\Ps_j\.\Jc\.\Ps_j\.H_j\Big\}
     =\innerP{\Rm_\la(F)}{H},
  \end{align*}
 where we used the equality $G_{kj}(\bla)=G^*_{jk}(\la)$ valid for all $j\neq k$ and $G_{jj}(\bla)=G_{jj}^*(\la)-\Jc$ 
 for any $j\in\oinftyZ$, see \cite[Lemma~4.2]{slC.pZ15}. Therefore, equality~\eqref{E:P.A3.3} reduces 
 to~\eqref{E:P.A3.1} immediately.
\end{proof}

In the remaining part of this section we establish several properties of $M_+(\la)$-function which will be crucial in 
the next section.

\begin{lemma}\label{L:MR1B}
 Let Hypothesis~\ref{H:basic} be satisfied. If $M_+(\la)$ is holomorphic on $(a,b)\subseteq\Rbb$, then
  \begin{enumerate}[leftmargin=10mm,topsep=1mm,label={{\normalfont{(\roman*)}}}]
%    \item $\Xc^+(\la)$ is holomorphic in $(a,b)\subseteq\Rbb$;
   \item $M_+'(\la)>0$ for all $\la\in(a,b)$;
   \item $\Xc^+(\la)\in\ltp$ with $\normP{\Xc^+(\la)}\leq M'(\la)$ for all $\la\in(a,b)$.
  \end{enumerate}
\end{lemma}
\begin{proof}
 It suffices to show only the second part, because then the inequality $M_+'(\la)>0$ is a simple consequence of the 
 strong Atkinson condition. Let $\mu\in\Cbb\stm\Rbb$ be arbitrary. By Theorem~\ref{T:Lagrange} we have 
  \begin{equation*}%\label{E:MR4}
   (\bar{\mu}-\mu)\.\sum_{k=0}^\infty \Xc^{+*}_k(\mu)\.\Ps_k\.\Xc^+(\mu)=M_+^*(\mu)-M_+(\mu),
  \end{equation*}
 i.e., 
  \begin{equation*}
   \normP{\Xc^+(\mu)}^2=\frac{M_+(\mu)-M_+^*(\mu)}{2i\.\im(\mu)}.
  \end{equation*}
 Hence for a fixed $\la\in(a,b)$ and $\nu>0$ small enough we obtain
  \begin{equation*}
   \normP{\Xc^+(\la+i\.\nu)}^2=\frac{M_+(\la+i\.\nu)-M_+^*(\la+i\.\nu)}{2i\.\nu}
  \end{equation*}
 with the right-hand side tending to $M_+'(\la)$ as $\nu\to0^+$, because $M_+(\la)$ is holomorphic, in which case the 
 symmetric derivative coincides with the classical derivative. Consequently,
  \begin{equation*}
   \sum_{k=0}^N \Xc^{+*}_k(\la)\.\Ps_k\.\Xc^+(\la)\leq M_+'(\la)
  \end{equation*}
 for any $N\in\oinftyZ$, which implies
  \begin{equation*}
   0<\normP{\Xc^+(\la)}^2\leq M_+'(\la).\qedhere
  \end{equation*}
\end{proof}

Now let us recall the definition of the {\it norm} of the resolvent relation given as 
 \begin{equation*}
  \norm{\Rc_{\TLP}(\la)}\coloneq 
  \inf\big\{C\geq0\mid \normP{f}\leq C\.\normP{z}\ \text{ for all } \{\zclass,\fclass\}\in\Rc_{\TLP}(\la)\big\}
 \end{equation*}
for every $\la\in\rh(\TLP)$, see \cite[Section~2.7]{pZ.slC:SAE2}.

\begin{lemma}\label{L:MR1A} 
 Let Hypothesis~\ref{H:basic} be satisfied. For any $f,g\in\ltp$ and $\la\in\rh(\TLP)$ the function 
 $h(\la)\coloneq \innerP{\gclass}{\Rc_{\TLP}(\la)(f)}$ exists and it is holomorphic with respect to $\la$.
\end{lemma}
\begin{proof}
 Let $\la\in\rh(\TLP)$ be arbitrary. Since $\rh(\TLP)$ is an open set, there is a neighborhood $\Oc(\la)$ with the 
 radius less than $\norm{\Rc_{\TLP}(\la)}^{-1}$ such that $\Oc(\la)\subseteq\rh(\TLP)$. Then for any $\mu\in\Oc(\la)$ 
 we have
  \begin{align}
   h(\mu)
    &=\innerP{\gclass}{\Rc_{\TLP}(\mu)(f)}
     =\sum_{k=0}^\infty g_k^*\.\Ps_k\.\Rc_{\TLP}(\mu)(f)_k \notag\\
    &=\sum_{k=0}^\infty g_k^*\.\Ps_k\bigg[\sum_{s=0}^\infty (\mu-\la)^s\.\Rc_{\TLP}^{s+1}(\la)(f)_k\bigg]
     =\sum_{k=0}^\infty\sum_{s=0}^\infty (\mu-\la)^s\.g_k^*\.\Ps_k\.\Rc_{\TLP}^{s+1}(\la)(f)_k,\label{E:L.MR1A.1}
   \end{align}
 with an arbitrary representative of $\Rc_{\TLP}(\la)(f)$ in the sums and its Neumann series representation, see 
 \cite[Theorem~2.44]{pZ.slC:SAE2}. Since $\abs{\mu-\la}\times \norm{\Rc_{\TLP}(\la)}<1$ and
  \begin{align*}
   &\sum_{k=0}^\infty\sum_{s=0}^\infty \abs[\big]{(\mu-\la)^s\.g_k^*\.\Ps_k\.\Rc_{\TLP}^{s+1}(\la)(f)_k}
    \leq \sum_{s=0}^\infty \abs{\mu-\la}^s\times\normP{g}\times \normP{\Rc_{\TLP}^{s+1}(\la)(f)}\\
    &\hspace*{16mm}\leq \normP{\Rc_{\TLP}(\la)(f)}\times\normP{g}\times
      \sum_{s=0}^\infty \abs{\mu-\la}^s\times \norm{\Rc_{\TLP}(\la)}^s\. \normP{f}<\infty,
  \end{align*}
 the existence of $h(\mu)$ follows. Consequently, we can also interchange the order in the double sum 
 in~\eqref{E:L.MR1A.1}, which leads to the expression
   \begin{align*} 
    h(\mu)=\sum_{s=0}^\infty\sum_{k=0}^\infty (\mu-\la)^s\.g_k^*\.\Ps_k\.\Rc_{\TLP}^{s+1}(\la)(f)_k
          =\sum_{s=0}^\infty (\mu-\la)^s\.\innerP{\gclass}{\Rc_{\TLP}^{s+1}(\la)(f)},
  \end{align*}
 i.e., $h(\mu)$ is represented be a convergent power series centered at $\la$, which guarantees its holomorphic 
 property with respect to $\la$. 
\end{proof}

In~\cite[Theorem~3.7]{pZ21} we proved that for any finite discrete interval $\onZ$, matrix $\be\in\Ga$, and 
number $\la\in\Cbb\stm\Rbb$, the regular Weyl--Titchmarsh function $M_{N+1}(\la)$ admits the Riemann--Stieltjes 
integral representation
 \begin{align*}%\label{E:MR8}
  M_{N+1}(\la)&=-[\be\.\tZ_{N+1}(\la)]^{-1}\.\be\.\hZ_{N+1}(\la)
              =M^{[0]}+\la\.M^{[1]}
                      +\int_{-\infty}^\infty\bigg(\frac{1}{t-\la}-\frac{t}{1+t^2}\bigg)\.\dtau_{\al,\be}(t)
 \end{align*}
with
 \begin{equation*}
  M^{[0]}\coloneq \re M_{N+1}(i) \qtextq{and} M^{[1]}\coloneq \lim_{\mu\to\infty} M_{N+1}(i\mu)/(i\mu)\geq0
 \end{equation*}
and $\tau_{\al,\be}$ being a right continuous $n\times n$ matrix-valued step function with jumps occurring 
only at the eigenvalues of~\eqref{E:Sla} as well as possessing the Hermitian property
 \begin{equation*}
  \tau^*_{\al,\be}(t)=\tau_{\al,\be}(t) \qtext{for all $t\in\Rbb$}
 \end{equation*}
and the nondecreasing property
 \begin{equation*}
  \tau_{\al,\be}(t)\leq \tau_{\al,\be}(s) \qtext{for any $t\leq s$.}
 \end{equation*}
Furthermore, we have also
 \begin{equation}\label{E:Im.Mla.integral}
  \frac{\im M_{N+1}(\la)}{\im(\la)}=M^{[1]}+\int_{-\infty}^\infty \frac{1}{\abs{t-\la}^2}\.\dtau_{\al,\be}(t)
 \end{equation}
and the nested property of the Weyl disks yields
 \begin{equation*}
  0\leq \int_{-\infty}^\infty (1+t^2)^{-1} \.\dtau_{\al,\be}(t)\leq T_0
 \end{equation*}
for some $T_0\in\Cbb^{n\times n}$ being independent of $N$, see \cite[pp.~33--34]{pZ21}. Consequently, there is a 
nondecreasing right continuous $n\times n$ Hermitian matrix-valued {\it limiting spectral function} $\tau$ and it holds
 \begin{equation}\label{E:MR9AA}
  \int_{-\infty}^\infty (1+t^2)^{-1} \.\dtau(t)\leq T_0
 \end{equation}
by the Helly's convergence theorem, see~\cite[Theorem~5.7.6]{gaM.aS.mT19}. The following statement shows that similar 
integral representation is valid also for $M_+(\la)$.

\begin{lemma}\label{L:MR2alpha}
 Let Hypothesis~\ref{H:basic} be satisfied. For any $\la\in\Cbb\stm\Rbb$ we have
  \begin{equation}\label{E:MR9A}
   M_+(\la)=M_+^{[0]}+\la\.M_+^{[1]}+\int_{-\infty}^\infty \Big(\frac{1}{t-\la}-\frac{t}{1+t^2}\Big)\.\dtau(t)
  \end{equation}
 with $M_+^{[0]}\coloneq \re M_+(i)$ and 
 $M_+^{[1]}\coloneq \im M_+(i)-\int_{-\infty}^\infty \frac{1}{1+t^2}\.\dtau(t)$.
\end{lemma}
\begin{proof}
 Let $R>0$ and $N\in\oinftyZ$ be arbitrary. Then
  \begin{align*}
   \Ups
    &\coloneq\frac{\im M_+(\la)}{\im(\la)}-\im M_+(i)
      -\int_{\abs{t}\leq R} \bigg(\frac{1}{\abs{t-\la}^2}-\frac{1}{1+t^2}\bigg)\.\dtau(t)\\
    &\hspace*{0.9mm}
     =-\int_{\abs{t}\leq R} 
               \bigg(\frac{1}{\abs{t-\la}^2}-\frac{1}{1+t^2}\bigg)\.\big(\dtau(t)-\dtau_{\al,\be}(t)\big)\\
     &\hspace*{15mm}+\int_{\abs{t}>R} 
               \bigg(\frac{1}{\abs{t-\la}^2}-\frac{1}{1+t^2}\bigg)\.\dtau_{\al,\be}(t)\\
      &\hspace*{27mm}+\bigg[\frac{\im M_+(\la)-\im M_{N+1}(\la)}{\im\la}-\im M_+(i)+\im M_{N+1}(i)\bigg]
  \end{align*}
 by~\eqref{E:Im.Mla.integral}. The first and third integrals on the right-hand side tend to $0$ as $N\to\infty$ because 
 of the uniformly bounded variation of $\tau_{\al,\be}$ on $[-R,R]$, the definition of $M_+(\la)$, and the limit point 
 case hypothesis. If we put $\la\coloneq \la_0+i\.\nu$, then it also holds
  \begin{equation*}
   \abs[\bigg]{\frac{1}{\abs{t-\la}^2}-\frac{1}{1+t^2}}
   =\abs[\bigg]{\frac{1+2\.\la_0\.t-\la_0^2-\nu^2}{(1+t^2)\.((t-\la_0)^2+\nu^2)}}
   \leq \frac{1}{1+t^2}\.\frac{\ga}{\abs{t}}
  \end{equation*}
 for all $t\in\Rbb$ and a suitable $\ga>0$, because $\abs{t\.(1+2\.\la_0\.t-\la_0^2-\nu^2)}/((t-\la_0)^2+\nu^2)$ is 
 continuous on $\Rbb$ and
  \begin{equation*}
   \lim_{t\to\pm\infty} \frac{\abs{t(1+2\.\la_0\.t-\la_0^2-\nu^2)}}{((t-\la_0)^2+\nu^2)}=2\.\abs{\la_0}.
  \end{equation*}
 So,
  \begin{equation*}
   \int_{\abs{t}>R} \bigg(\frac{1}{\abs{t-\la}^2}-\frac{1}{1+t^2}\bigg)\.\dtau_{\al,\be}(t)
    \leq \int_{\abs{t}>R} \frac{1}{1+t^2}\.\frac{\ga}{\abs{t}}\.\dtau_{\al,\be}(t)
    \leq \frac{\ga}{R}\.T_0,
  \end{equation*}
 which does not depend on $N$. If we fix $R>0$, then for $N\to\infty$ we get 
  \begin{equation*}
   \abs{\Ups}\leq \frac{\ga}{R}\.T_0,
  \end{equation*}
 and thus $\abs{\Ups}\to0$ as $R\to\infty$, which shows that
  \begin{equation*}
   \im M_+(\la)=\im(\la)\.\bigg(\im M_+(i)-\int_{-\infty}^\infty \frac{1}{1+t^2}\.\dtau(t)\bigg)
                 +\int_{-\infty}^\infty \frac{\im(\la)}{\abs{t-\la}^2}\.\dtau(t).
  \end{equation*}
 Consequently, the holomorphic property of $M_+(\la)$ on $\Cbb\stm\Rbb$ implies that it is uniquely determined by its 
 imaginary part up to an additive Hermitian matrix as in the proof of~\cite[Theorem~3.7]{pZ21}, which justifies
 equality~\eqref{E:MR9A}.
\end{proof}

As the next result, we derive a certain limiting property of $M_+(\la)$.

\begin{lemma}\label{L:MR2A}
 Let Hypothesis~\ref{H:basic} be satisfied. For any $\la_0\in\Rbb$ we have
  \begin{equation}\label{E:MR7}
   \lim_{\nu\to0} \nu\.M_+(\la_0+i\.\nu)=i\.\big(\tau(\la_0)-\tau(\la_0^-)\big).
  \end{equation}
\end{lemma}
\begin{proof}
 From~\eqref{E:MR9A} we have
  \begin{equation}\label{E:MR10}
   \begin{aligned}
    \nu\.M(\la_0+i\.\nu)
     &=\nu\.M_+^{[0]}+\nu\.(\la_0+i\.\nu)\.M_+^{[1]} \\
       &\hspace*{15mm}+\nu\int_{-\infty}^\infty \Big(\frac{t-\la_0}{(t-\la_0)^2+\nu^2}-\frac{t}{1+t^2}\Big)\.\dtau(t)\\
         &\hspace*{30mm}+i\.\int_{-\infty}^\infty \frac{\nu^2}{(t-\la_0)^2+\nu^2}\.\dtau(t)
   \end{aligned}
  \end{equation}
 for any $\la_0\in\Rbb$ and $\nu\in\Rbb\stm\{0\}$. Let us fix $\de>0$ and put $\tau(t)\coloneq\tau_1(t)+\tau_2(t)$ with 
 the matrix-valued step function
  \begin{equation*}
   \tau_1(t)\coloneq\begin{cases}
                     \tau(\la_0),& t\geq\la_0,\\[1mm]
                     \tau(\la_0^-),& t<\la_0,
                    \end{cases}
  \end{equation*}
 possessing the only jump $\tau(\la_0)-\tau(\la_0^-)$ at $t=\la_0$ and with the second component 
 $\tau_2(t)\coloneq \tau(t)-\tau_1(t)$ being continuous at $t=\la_0$. We start with the last integral, which we split 
 into two parts as
  \begin{align*}%\label{E:MR11}
   \int_{-\infty}^\infty \frac{\nu^2}{(t-\la_0)^2+\nu^2}\.\dtau(t)
    &=\int_{\la_0-\de}^{\la_0+\de} \frac{\nu^2}{(t-\la_0)^2+\nu^2}\.\dtau(t)
      +\int_{\abs{t-\la_0}\geq\de} \frac{\nu^2}{(t-\la_0)^2+\nu^2}\.\dtau(t).
  \end{align*}
 Then, by the definition of the Riemann--Stieltjes integral, we have
  \begin{equation}\label{E:MR12}
   \nu^2\.\int_{\abs{t-\la_0}\leq\de} 
    \frac{1}{(t-\la_0)^2+\nu^2}\.\dtau_1(t)=\tau(\la_0)-\tau(\la_0^-).
  \end{equation}
 Simultaneously, since
  \begin{equation*}
   \frac{\nu^2}{(t-\la_0)^2+\nu^2}\leq1
  \end{equation*}
 for all $t\in[\la_0-\de,\la_0+\de]$ and any $\nu\in\Rbb$, it follows from the Osgood's convergence 
 theorem, see \cite[Theorem~5.7.5]{gaM.aS.mT19} or \cite[Corollary on p.~274]{jdW51}, that
  \begin{equation*}
   \lim_{\nu\to0} \int_{\la_0-\de}^{\la_0+\de} \frac{\nu^2}{(t-\la_0)^2+\nu^2}\.\dtau_2(t)
    =\int_{\la_0-\de}^{\la_0+\de} \lim_{\nu\to0} \frac{\nu^2}{(t-\la_0)^2+\nu^2}\.\dtau_2(t)=0.
  \end{equation*}
 It remains to evaluate the integral
  \begin{equation*}
   \lim_{\nu\to0} \nu^2\int_{\abs{t-\la_0}\geq\de} \frac{1}{(t-\la_0)^2+\nu^2}\.\dtau(t).
  \end{equation*}
 Since $\abs{t-\la_0}\geq\de$ and we have the upper bound
  \begin{equation*}
   \frac{1+t^2}{(t-\la_0)^2+\nu^2}
    \leq \frac{1}{(t-\la_0)^2}+\Big(\frac{t}{t-\la_0}\Big)^2
    \leq \frac{1+(\abs{\la_0}+\de)^2}{\de^2}\eqcolon\ga,
  \end{equation*} 
 we obtain the estimate
  \begin{equation*}%\label{E:MR14}
   \frac{1}{(t-\la_0)^2+\nu^2}\leq\frac{\ga}{1+t^2},
  \end{equation*}
 which implies
  \begin{align*}
   0\leq \lim_{\nu\to0} \nu^2 \int_{\abs{t-\la_0}\geq\de} \frac{1}{(t-\la_0)^2+\nu^2}\.\dtau(t)
    &\leq \ga\.\lim_{\nu\to0} \nu^2 \int_{-\infty}^\infty \frac{1}{1+t^2}\.\dtau(t)
     =\ga\.\lim_{\nu\to0} \nu^2\.T_0=0
  \end{align*}
 by~\eqref{E:MR9AA}. Finally, we show that also
  \begin{equation*}
   \lim_{\nu\to0} \nu\.\int_{-\infty}^\infty \Big(\frac{t-\la_0}{(t-\la_0)^2+\nu^2}-\frac{t}{1+t^2}\Big)\.\dtau(t)=0,
  \end{equation*}
 for which we may proceed in a similar way as above. More precisely,
  \begin{equation*}
   \int_{\la_0-\de}^{\la_0+\de}\frac{t-\la_0}{(t-\la_0)^2+\nu^2}\.\dtau_1(t)=0
  \end{equation*}
 by definition. Since $\lim_{\nu\to0} \nu\.(t-\la_0)/[(t-\la_0)^2+\nu^2]=0$ and $\nu\.(t-\la_0)/[(t-\la_0)^2+\nu^2]$ 
 is bounded on $[\la_0-\de,\la_0+\de]$ by the inequality of arithmetic and geometric means, it follows again from the 
 Osgood's convergence theorem that 
  \begin{equation*}
   \lim_{\nu\to0} \int_{\la_0-\de}^{\la_0+\de}\frac{\nu\.(t-\la_0)}{(t-\la_0)^2+\nu^2}\.\dtau_2(t)=0.
  \end{equation*}
 Next, by a direct calculation
  \begin{equation*}
   \lim_{\nu\to0} \nu\.\int_{\la_0-\de}^{\la_0+\de}\frac{t}{1+t^2}\.\dtau(t)
   \leq \lim_{\nu\to0}\nu\.(\la_0+\de)\.\int_{\la_0-\de}^{\la_0+\de}\frac{1}{1+t^2}\.\dtau(t)=0
  \end{equation*}
 and simultaneously
  \begin{equation}\label{E:MR15}
   \lim_{\nu\to0} \nu\.\int_{\abs{t-\la_0}\geq\de} 
                    \bigg(\frac{t-\la_0}{(t-\la_0)^2+\nu^2}-\frac{t}{1+t^2}\bigg)\.\dtau(t)=0,
  \end{equation}
 because for $\abs{\nu}\leq1$ we get
  \begin{equation*}
   \abs[\bigg]{\frac{\nu\.(t-\la_0)}{(t-\la_0)^2+\nu^2}-\frac{\nu\.t}{1+t^2}}
    \leq \abs[\bigg]{\frac{t-\la_0}{(t-\la_0)^2+\nu^2}}+\abs[\bigg]{\frac{\nu\.t}{1+t^2}}
    \leq (1+\abs{\nu})/2\leq1
  \end{equation*}
 and for any $t\in\Rbb$ it holds
  \begin{equation*}
   \lim_{\nu\to0} \nu\bigg(\frac{t-\la_0}{(t-\la_0)^2+\nu^2}-\frac{t}{1+t^2}\bigg)=0,
  \end{equation*}
 which together with the Osgood's convergence theorem yields
  \begin{equation*}
   \lim_{\nu\to0} \nu\.\int_{M\geq\abs{t-\la_0}\geq\de} 
                    \bigg(\frac{t-\la_0}{(t-\la_0)^2+\nu^2}-\frac{t}{1+t^2}\bigg)\.\dtau(t)=0.
  \end{equation*}
 As the latter equality is valid for any $M>\de$, equality~\eqref{E:MR15} is proven. Therefore, the only non-zero term 
 in~\eqref{E:MR10} as $\nu\to0$ comes from~\eqref{E:MR12}, which leads to~\eqref{E:MR7}.
\end{proof}

For any $N\in\oinftyZ$ and $\be\in\Ga$ we already know that the function $M_{N+1}(\la)$ is holomorphic on $\Cbb_+$ and 
$\Cbb_-$ and, in addition, it is even a {\it Nevanlinna} (or {\it Herglotz}) {\it function}, see 
\cite[Section~5]{fG.erT00}. Thus, it has only isolated singularities being simple poles located on the real line. 
Actually, these singularities are the eigenvalues of the associated eigenvalue problem. On the other hand, if
$\la\in\Rbb$ is not an eigenvalue or equivalently $M_{N+1}(\la)$ is holomorphic at $\la$, then the spectral function 
$\tau_{\al,\be}(t)$ is constant in a neighborhood $\Oc(\la)\subseteq\Rbb$. This is a general property valid for every 
Nevanlinna function by \cite[Proposition~A.4.5]{jB.sH.hsvdS20}, and so it is true also for the limiting 
Weyl--Titchmarsh function $M_+(\la)$ and the limiting spectral function $\tau(t)$, because $M_+(\la)$ is a Nevanlinna 
function as well. The latter fact follows from the equality $M_+^*(\la)=M_+(\bla)$ and inequality $\sgn(\im \la)\.\im 
M_+(\la)\geq0$ for all $\la\in\Cbb\stm\Rbb$ and it fully agrees with the explicit integral representation established 
in Lemma~\ref{L:MR2alpha}. Furthermore, as the last result of this section we show that $M_+(\la)$ may have isolated 
singularities of one type only, compare with \cite[Theorem~5.4(vi)]{fG.erT00} and \cite[p.~661]{jB.sH.hsvdS20}.

\begin{lemma}\label{L:MR2C}
 Let Hypothesis~\ref{H:basic} be satisfied. If $\la_0$ is an isolated singularity of $M_+(\la)$, then $\tau(\cdot)$ 
 has a nonzero jump at $\la_0$ and this point is a simple pole of $M_+(\la)$.
\end{lemma}
\begin{proof}
 If $\la_0$ is an isolated singularity of $M_+(\la)$, then $\la_0\in\Rbb$ by the holomorphic property of $M_+(\la)$ on
 $\Cbb_+\cup\Cbb_-$. Thus, $M_+(\la)$ is holomorphic for all 
 $\la\in\Oc^*_\de(\la_0)\coloneq\Oc_\de(\la_0)\stm\{\la_0\}$ with a suitable $\de>0$, so $\tau(\cdot)$
 is constant on $(\la_0-\de,\la_0)$ and $(\la_0,\la_0+\de)$ by the preceding discussion. Simultaneously, 
 $\tau(\cdot)$ cannot be continuous at $\la_0$, because otherwise $\tau(\cdot)$ would be constant on 
 $(\la_0-\de,\la_0+\de)$ and, consequently, $M_+(\la)$ would be holomorphic on $\Oc_\de(\la_0)$ by the same reason. 
 Therefore, $\tau(\cdot)$ has a~jump discontinuity at $\la_0$ and furthermore, for any $\la\in\Oc^*_\de(\la_0)$, we 
 have the expansion
  \begin{align*}
   M_+(\la)
   &=M_+^{[0]}+\la\.M_+^{[1]}
    +\bigg(\int_{-\infty}^{\la_0-\de}+\int_{\la_0+\de}^\infty\bigg)\bigg(\frac{1}{t-\la} 
                                                                            +\frac{t}{1+t^2}\bigg)\.\dtau(t)\\
    &\hspace*{60mm}+\int_{\la_0-\de}^{\la_0+\de}\bigg(\frac{1}{t-\la}+\frac{t}{1+t^2}\bigg)\.\dtau(t)\\
   &=M_+^{[0]}+\la\.M_+^{[1]}
    +\bigg(\int_{-\infty}^{\la_0-\de}+\int_{\la_0+\de}^\infty\bigg)\bigg(\frac{1}{t-\la}
                                                                            +\frac{t}{1+t^2}\bigg)\.\dtau(t)\\
    &\hspace*{43mm}-\frac{\la_0}{1+\la_0^2}\big(\tau(\la_0)-\tau(\la_0^-)\big)
    +\frac{\tau(\la_0)-\tau(\la_0^-)}{\la_0-\la},
  \end{align*}
 i.e.,
  \begin{align*}
   M_+(\la)-\frac{\tau(\la_0)-\tau(\la_0^-)}{\la_0-\la}
    &=M_+^{[0]}+\la\.M_+^{[1]}\\
      &\hspace*{15mm}+\bigg(\int_{-\infty}^{\la_0-\de}+\int_{\la_0+\de}^\infty\bigg)\bigg(\frac{1}{t-\la}
                                                                            +\frac{t}{1+t^2}\bigg)\.\dtau(t)\\
    &\hspace*{50mm}-\frac{\la_0}{1+\la_0^2}\big(\tau(\la_0)-\tau(\la_0^-)\big)\eqcolon g(\la)
  \end{align*}
 with $g(\la)$ being holomorphic on $\Oc^*_\de(\la_0)$. We show that $g(\la)$ is bounded at $\la_0$. Let 
 $[\la_1,\la_2]\subseteq\Oc_\de(\la_0)$ be an arbitrary interval with $\la_1<\la_0<\la_2$. Then for any 
 $\la\in[\la_1,\la_2]$ we have
  \begin{equation*}
   \int_{-\infty}^{\la_0-\de} \bigg(\frac{1}{t-\la}+\frac{t}{1+t^2}\bigg)\.\dtau(t)
   \leq \int_{-\infty}^{\la_0-\de} \abs[\Big]{\frac{t\.\la}{t-\la}}\.\frac{1}{1+t^2}\.\dtau(t)\leq \ga\.T_0
  \end{equation*}
 for a suitable $\ga>0$ and similarly we can show the boundedness of the second integral. Thus,
  \begin{equation*}
   \lim_{\la\to\la_0} (\la-\la_0)\.M_+(\la)
    =\lim_{\la\to\la_0} (\la-\la_0)\.g(\la)-\tau(\la_0)+\tau(\la_0^-)
    =-\tau(\la_0)+\tau(\la_0^-)\neq0,
  \end{equation*}
 which means that $\la_0$ is a simple pole of $M_+(\la)$.
\end{proof}

%%%%%%%%%%%%%%%%%%%%%%%%%%%%%%%%%%%%%%%%%%%% SECTION %%%%%%%%%%%%%%%%%%%%%%%%%%%%%%%%%%%%%%%%%%%%%%%%%%%%%%%%%%%%%%%%%

\section{Main results}\label{S:spectrum}

At this moment we have all tools, which we need for the proof of Theorem~\ref{T:intro.spectrum} stated in the 
introduction. For simplicity, the conclusions are proven separately.

\begin{theorem}\label{T:MR1}
 Let Hypothesis~\ref{H:basic} be satisfied. We have $\la_0\in\rh(\TLP)$ if and only if $M_+(\la)$ is holomorphic at 
 $\la_0$. Furthermore, the resolvent relation admits the representation
  \begin{equation}\label{E:MR1}
   \Rc_{\TLP}(\la_0)\coloneq(\TLP-\la_0\.I)^{-1}
    =\Bigg\{\bigg\{\fclass,\Big[\sum_{j=0}^\infty G_{kj}(\la_0)\.\Ps_j\.f_j\Big]\bigg\}\,\Big|\, f\in\ltp\Bigg\},
  \end{equation}
 where
  \begin{equation*}%\label{E:MR2}
   G_{kj}(\la_0)\coloneq
    \begin{cases}
     \tZ_k(\la_0)\.\Xc^{+*}_j(\bla_0),& k\in[0,j]_\sZbb,\\
     \Xc^+_k(\la_0)\.\tZ^*_j(\bla_0),& k\in[j+1,\infty)_\sZbb.
    \end{cases}
  \end{equation*}
\end{theorem}
\begin{proof}
 If $\la_0\in\Cbb\stm\Rbb\subseteq\rh(\TLP)$, then the first part follows from Montel's theorem. More precisely, for 
 any $\la\in\Cbb\stm\Rbb$ and $k\in\oinftyZ$ we have 
  \begin{equation*}
   M_k(\la,\al,\be)\in C_k(\la)\coloneq
   \big\{M\in\Cbb^{n\times n}\mid \Xc_k^*(\la,M)\,\Jc\Xc_k(\la,M)=0\big\}
  \end{equation*}
 by~\cite[Theorem~3.2]{rSH.pZ14:JDEA} with the (regular) Weyl solution  $\Xc(\la,M)\coloneq\hZ(\la)+\tZ(\la)\.M(\la)$, 
i.e., it is located on the Weyl circle being the boundary of the 
 Weyl disk. At the same time, $M_k(\cdot,\al,\be)$ is holomorphic on $\Cbb_+\cup\Cbb_-$ for any $k\in\oinftyZ$ and it 
 is pointwise convergent for any $\la\in\Cbb_+\cup\Cbb_-$ to the center of the limiting Weyl disk $P_+(\la)$, i.e., 
 $M_+(\la)\coloneq\lim_{k\to\infty} M_k(\la,\al,\be)=P_+(\la)$. Moreover, the nested property of Weyl disks 
 $D_k(\la)\coloneq \big\{M\in\Cbb^{n\times n}\mid \sgn \im(\la)\,i\,\Xc_k^*(\la,M)\,\Jc\Xc_k(\la,M)\leq0\big\}$ yields
  \begin{equation*}
   M_k(\la)\in C_k(\la)\subseteq D_k(\la)\subseteq D_{k-1}(\la)\subseteq\cdots\subseteq D_0(\la),\\
  \end{equation*}
 and, consequently, it can be also represented as
  \begin{equation*}
   M_k(\la)=P_{N_0+1}(\la)+R_{N_0+1}(\la)\.V_k\.R_{N_0+1}(\bla)
  \end{equation*}
 for any $k\in[N_0+1,\infty)_\sZbb$ and a suitable $V_k^* V_k\leq I$. Thereby, the sequence 
 $\{M_k(\la)\}\in\Cbb^{n\times n}(\oinftyZ)$ is locally uniformly bounded with
  \begin{align*}
   \normS{M_k(\la)}
    \leq\normS{P_{N_0+1}(\la)}+\normS{R_{N_0+1}(\la)}\times\normS{R_{N_0+1}(\bla)}\eqcolon \ga(\la)\leq \ga,
  \end{align*}
 and $\ga\coloneq \max\{\ga(\la)\mid \text{$\la$ belongs to a compact set}\}$. Therefore, $M_+(\la)$ is holomorphic on 
 $\Cbb_+\cup\Cbb_-$ and equality~\eqref{E:MR1} follows from~\cite[Theorem~4.3]{slC.pZ15}.
 
 Now, let $\la_0\in\rh(\TLP)\cap\Rbb$. Since $\rh(\TLP)$ is open, there is an interval $(a,b)$ such that 
 $\la_0\in(a,b)\subseteq\Rbb$. We show that $M_+(\la)$ can be analytically continued to $\la_0$, i.e., that 
 $M_+(\la_0)$ is holomorphic. Thus, let $\la\in\Cbb\stm\Rbb$ be arbitrary. Then we know that
  \begin{equation*}
   \lim_{k\to\infty} \Xc_k^{+*}(\la)\.\Jc\.\Xc_k^+(\nu)=0
  \end{equation*}
 for any $\nu\in\Cbb\stm\Rbb$ by~\cite[Equation~(4.3)]{slC.pZ15}, which together with the Lagrange identity yields
  \begin{equation}\label{E:MR3}
   (i-\la)\.\sum_{k=0}^\infty \Xc_k^{+*}(-i)\.\Ps_k\.\Xc_k^+(\la)
    =\Xc_0^{+*}\.\Jc\.\Xc^+_0(\la)
    =M_+(i)-M_+(\la).
  \end{equation}
 By definition, the columns of $\Xc^+(\la)$ and $\Xc^+(-i)$ belong to $\ltp$ and these two $2n\times n$ matrix-valued 
 solutions satisfy
  \begin{equation*}
   \al\.\Xc_0^+(\la)=I=\al\.\Xc_0^+(-i).
  \end{equation*}
 Thus, the matrix-valued function $Z\coloneq \Xc^+(\la)-\Xc^+(-i)$ is such that $\al\.Z_0=0$ and
  \begin{align*}
   Z_k=\Xc_k^+(\la)-\Xc_k^+(-i)
      &=\Sbb_k(\la)\.Z_{k+1}+(\la+i)\.\Vc_k\.\Xc_{k+1}^+(-i)\\
      &=\Sbb_k(\la)\.Z_{k+1}-(\la+i)\.\Jc\.\Ps_k\.\Xc_k^+(-i)
  \end{align*}
 for all $k\in\oinftyZ$ i.e., it solves system~\Slaf{\la}{f} with $f\coloneq (\la+i)\.\Xc^+(-i)$. Therefore, it can be 
 written as
  \begin{equation*}
   Z=(\la+i)\.\sum_{j=0}^\infty G_{\cdot j}(\la)\.\Ps_j\.\Xc_j^+(-i)=(\la+i)\.\Rm_\la(\Xc^+(-i))
  \end{equation*}
 and, consequently, 
  \begin{equation*}
   \Xc^+(\la)=\Xc^+(-i)+(\la+i)\.\Rm_\la(\Xc^+(-i)).
  \end{equation*}
 Upon inserting this expression into~\eqref{E:MR3} we obtain
  \begin{align*}
   M_+(i)-M_+(\la)
    &=(i-\la)\.\sum_{k=0}^\infty \Xc_k^{+*}(-i)\.\Ps_k\.\Xc_k^+(\la)\\
    &=(i-\la)\.\sum_{k=0}^\infty \Xc_k^{+*}(-i)\.\Ps_k\.\big[\Xc_k^+(-i)+(\la+i)\.\Rm_\la(\Xc^+(-i))_k\big]\\
    &=(i-\la)\.\normP{\Xc^+(-i)}^2
      -(1+\la^2)\.\innerP[\big]{[\Xc^+(-i)]}{\Rc_{\TLP}(\la)(\Xc^+(-i))},
  \end{align*}
 which together with Lemma~\ref{L:MR1A} implies the holomorphic property of $M_+(\la)$ at $\la_0$, i.e., the existence 
 of an analytic continuation into $\la_0$. 
 
 On the other hand, let $\la_0\in\Cbb$ be such that $M_+(\la)$ is holomorphic at $\la_0$. We already know that this is 
 true when $\la_0\in\Cbb\stm\Rbb$. So, let us assume that this is true also for some $\la_0\in\Rbb$. Then, 
 $\Xc^+(\la_0)\in\ltp$ by Lemma~\ref{L:MR1B} and we need to show that $\la_0\in\rh(\TLP)$, i.e.,
 $\ker(\TLP-\la_0\.I)=\{0\}$ and $\ran(\TLP-\la_0\.I)=\ltp$. In fact, it suffices to prove $\ran(\TLP-\la_0\.I)=\ltp$,
 because in that case we get from the self-adjointness of $\TLP$ that
  \begin{equation*}
   \{0\}=\big(\ran(\TLP-\la_0\.I)\big)^\bot
        =\ker(\TLP-\la_0\.I)^*=\ker(\TLP-\bla_0\.I)=\ker(\TLP-\la_0\.I).
  \end{equation*}
 Let $f\in\ltp$ be arbitrary and put $\hz\coloneq\sum_{j=0}^\infty G_{\cdot j}(\la_0)\.\Ps_j\.f_j$. Then $\hz$ 
 solves~\Slaf{\la_0}{f} and $\al\.\hz_0=0$, so it remains to verify that $\hz\in\ltp$, which will imply that  
 $[\hz]\in\dom\TLP$, i.e., $\big\{[\hz],\fclass\big\}\in\TLP-\la_0\.I$ or equivalently $\fclass\in\ran(\TLP-\la_0\.I)$. 
 First of all, the holomorphic property of $M_+(\la)$ at $\la_0$ implies the continuity of $M_+(\la)$ on some 
 neighborhood $\Oc(\la_0)$, so it holds
  \begin{equation*}
   \lim_{\nu\to0^+} M_+(\la_0+i\.\nu)=M_+(\la_0) \qtextq{as well as}
   \lim_{\nu\to0^+} M^*_+(\la_0-i\.\nu)=M^*_+(\la_0).
  \end{equation*}
 Since we also know that $M_+(\la_0+i\.\nu)=M_+^*(\overline{\la_0+i\.\nu})=M_+^*(\la_0-i\.\nu)$ for any $\nu>0$, it 
 follows that $M_+(\la_0)=M_+^*(\la_0)$. Consequently,
  \begin{equation*}
   \Xc^{+*}_k(\la_0)\.\Jc\.\Xc_k^+(\la_0)\equiv \Xc^{+*}_0(\la_0)\.\Jc\.\Xc_0^+(\la_0)
    =M_+(\la_0)-M_+^*(\la_0)=0
  \end{equation*}
 for any $k\in\oinftyZ$ by Theorem~\ref{T:Lagrange}. Then the square summability of $\hz$ can be shown by the 
 same arguments as in~\cite[pp.~793--795]{slC.pZ15}, in which we use the latter equality instead of the fact 
 $M_+(\la)\in D_+(\la)$ applied for $\la\in\Cbb\stm\Rbb$. More precisely, if we replace $f$ by 
 $f^{[r]}$ with $f^{[r]}_k\coloneq f_k$ for $k\in[0,r]_\sZbb$ and $f^{[r]}_k\coloneq 0$ for $k\in[r+1,\infty)_\sZbb$, 
 then the corresponding solution $\hz^{[r]}$ satisfies $\normP{\hz^{[r]}}\leq \normP{f}/(\im \la)$. As 
 $\hz^{[r]}$ converges uniformly to $\hz$ for $r\to\infty$, this leads to the estimate 
 $\normP{\hz}\leq \normP{f}/(\im \la)$, i.e., $\hz\in\ltp$.
\end{proof} 

Now, we focus on the spectrum of $\TLP$ and start with its isolated points.

\begin{theorem}\label{T:MR2}
 Let Hypothesis~\ref{H:basic} be satisfied. A number $\la_0\in\Rbb$ belongs to $\sip(\TLP)$, i.e., it is an isolated 
 eigenvalue of $\TLP$, if and only if $M_+(\la)$ has a simple pole at $\la_0$, which is also equivalent to the 
 existence of the Laurent expansion of $M_+(\la)$ in a neighborhood of $\la_0$ as
  \begin{equation}\label{E:MR5}
   M_+(\la)=K_{-1}\.(\la-\la_0)^{-1}+K_0+K_1\.(\la-\la_0)+\cdots,
  \end{equation}
 where $K_{-1},K_0,\dots$ are $n\times n$ Hermitian matrices with
  \begin{equation}\label{E:MR6}
   K_{-1}=-\big(\tau(\la_0)-\tau(\la_0^-)\big)\lneqq0.
  \end{equation}
 In this case, all columns of $\tZ(\la_0)\.K_{-1}$ belong to $\ltp$ and the nonzero columns are eigenfunctions of 
 $\TLP$ related to $\la_0$. Furthermore, the columns of 
  \begin{equation*}
   \hZ(\la_0)+\tZ(\la_0)\.K_0+\Big[\frac{\d}{\d\la} \tZ(\la)\Big]_{\la=\la_0}\.K_{-1}
  \end{equation*}
 belong also to $\ltp$. 
\end{theorem}
\begin{proof}
 By the definition of $\sip(\TLP)$ and Theorem~\ref{T:MR1}, a number $\la_0$ is an isolated singularity of $M_+(\la)$ 
 if any only if it is an isolated point of $\si(\TLP)=\Cbb\stm\rh(\TLP)\subseteq\Rbb$, i.e., it belongs to 
 $\sip(\TLP)$. In this case, the Laurent series expansion~\eqref{E:MR5} exists in a neighborhood of $\la_0$ with 
 $K_{-1}$ given in~\eqref{E:MR6} by Lemma~\ref{L:MR2C} and its proof. The Hermitian property of matrices $K_j$ is a 
 simple consequence of the equality $M_+^*(\la)=M_+(\bla)$. 
 
 Next we show that the columns of $\tZ(\la_0)\.K_{-1}$ belong to $\ltp$. By the Lagrange identity we have for any 
 $\nu\in\Rbb$ that
  \begin{equation}\label{E:MR15.i}
   \nu^2\.\sum_{k=0}^\infty \Xc_k^{+*}(\la_0+i\.\nu)\.\Ps_k\.\Xc_k^+(\la_0+i\.\nu)
    =\nu\.\im M_+(\la_0+i\.\nu)=-K_{-1}+\nu^2\.K_1+\dots
  \end{equation}
 The right-hand side of~\eqref{E:MR15.i} tends to $-K_{-1}$ as $\nu\to0$, while the left-hand side tends to 
 $\normP{\tZ(\la_0)\.K_{-1}}^2$ showing that $\normP{\tZ(\la_0)\.K_{-1}}^2\leq K_{-1}$, i.e., the columns of 
 $\tZ(\la_0)\.K_{-1}$ truly belong to $\ltp$. Indeed, for any $m\in\oinftyZ$ we have
  \begin{equation}\label{E:MR16}
   \nu^2\.\sum_{k=0}^m \Xc_k^{+*}(\la_0+i\.\nu)\.\Ps_k\.\Xc_k^{+}(\la_0+i\.\nu)\leq -K_{-1}+\O(\nu^2)
  \end{equation}
 and simultaneously
  \begin{equation}\label{E:MR17}
   \begin{aligned}
    &\nu\.\Xc_k^{+}(\la_0+i\.\nu)\overset{\eqref{E:MR5}}{=}
     \nu\.\hZ_k(\la_0+i\.\nu)-i\.\tZ_k(\la_0+i\.\nu)\.K_{-1}\\
     &\hspace*{50mm}+\tZ_k(\la_0+i\.\nu)\.\big(\nu\.K_0+i\.\nu^2\.K_1-\nu^3\.K_2+\dots\big),
   \end{aligned}  
  \end{equation}
 which tends to $-i\.\tZ_k(\la_0)\.K_{-1}$ for all $k\in\oinftyZ$ as $\nu\to0$. Thus,
  \begin{align*}
   &\lim_{\nu\to0} \nu^2\.\sum_{k=0}^m \Xc_k^{+*}(\la_0+i\.\nu)\.\Ps_k\.\Xc_k^{+}(\la_0+i\.\nu)\\
    &\hspace*{40mm}=\sum_{k=0}^m \big(\tZ_k(\la_0)\.K_{-1}\big)^*\.\Ps_k\.\big(\tZ_k(\la_0)\.K_{-1}\big)\leq K_{-1}
  \end{align*}
 by~\eqref{E:MR16}, and so this uniform boundedness implies that
  \begin{equation}\label{E:MR18}
   \normP{\tZ(\la_0)\.K_{-1}}^2
    =\sum_{k=0}^\infty \big(\tZ_k(\la_0)\.K_{-1}\big)^*\.\Ps_k\.\big(\tZ_k(\la_0)\.K_{-1}\big)
    \leq K_{-1}
  \end{equation}
 as claimed above. Since also $\al\.\tZ_0(\la_0)\.K_{-1}=0$, it follows that all the nonzero columns of 
 $\tZ(\la_0)\.K_{-1}$ are eigenfunctions of $\TLP$ corresponding to $\la_0$.
 
 Finally, we prove the square summability of the columns of 
  \begin{equation*}
   \hZ(\la_0)+\tZ(\la_0)\.K_0+\Big[\frac{\d}{\d\la} \tZ(\la)\Big]_{\la=\la_0}\.K_{-1}.
  \end{equation*}
 For $\nu>0$, let us define
  \begin{equation*}
   Y(\nu)\coloneq \Xc^+(\la_0+i\.\nu)-(i\.\nu)^{-1}\.\tZ(\la_0)\.K_{-1}.
  \end{equation*}
 Then the columns of $Y(\nu)$ belong to $\ltp$ by the previous part and the square summability of the columns of the 
 Weyl solution $\Xc^+(\la_0+i\.\nu)$. Furthermore, by a direct calculation we can verify that
  \begin{equation*}
   \mL(Y(\nu))_k=(\la_0+i\.\nu)\.\Ps_k\.Y_k(\nu)+\Ps_k\.\tZ_k(\la_0)\.K_{-1},
  \end{equation*}
 i.e., $Y(\nu)$ solves the nonhomogeneous system~\Slaf{\la_0+i\.\nu}{f} with $f\coloneq \tZ(\la_0)\.K_{-1}$ and the 
 extended Lagrange identity yields
  \begin{equation}\label{E:MR17A}
   \begin{aligned}
    \!\!\!
    &Y^*_k(\nu)\.\Jc\.Y_k(\nu)\Big|_{0}^\infty\\
     &\hspace*{10mm}=-2i\.\sum_{k=0}^\infty \Big[\nu\.Y_k^*(\nu)\.\Ps_k\.Y_k(\nu)\\ 
      &\hspace*{28mm}+\im\Big(\big(\Xc^{+}_k(\la_0+i\.\nu)-(i\.\nu)^{-1}\tZ_k(\la_0)\.K_{-1}\big)^*\times
                                 \Ps_k\.\tZ_k(\la_0)\.K_{-1}\Big)\Big].
   \end{aligned}
  \end{equation}
 Upon applying the Lagrange identity once again, we obtain
  \begin{align*}
   \sum_{k=0}^\infty \Xc_k^{+*}(\la_0+i\.\nu)\.\Ps_k\.\tZ_k(\la_0)\.K_{-1}
    &=-(i\.\nu)^{-1}\.\Xc^{+*}_k(\la_0+i\.\nu)\.\Jc\.\tZ_k(\la_0)\.K_{-1}\Big|_{0}^\infty\\
    &=(i\.\nu)^{-1}\.\Xc^{+*}_0(\la_0+i\.\nu)\.\Jc\.\tZ_0(\la_0)\.K_{-1}
     =(i\.\nu)^{-1}\.K_{-1}
  \end{align*}
 or equivalently
  \begin{equation*}%\label{E:MR20}
   \nu\.\sum_{k=0}^\infty \Xc_k^{+*}(\la_0+i\.\nu)\.\Ps_k\.\tZ_k(\la_0)\.K_{-1}=-i\.K_{-1}.
  \end{equation*}
 Therefore, by letting $\nu\to0^+$ we get
  \begin{equation*}%\label{E:MR21}
   -i\.K_{-1}=\lim_{\nu\to0^+} \nu\.\sum_{k=0}^\infty \Xc_k^{+*}(\la_0+i\.\nu)\.\Ps_k\.\tZ_k(\la_0)\.K_{-1}
    =i\.\normP{\tZ_k(\la_0)\.K_{-1}}^2
  \end{equation*}
 by~\eqref{E:MR17} and the Moore--Osgood theorem, see \cite[Theorem~7.11]{wR76}. This means that in~\eqref{E:MR18} the 
 equality occurs. Consequently, equality~\eqref{E:MR17A} reduces to
  \begin{align*}
   Y^*_k(\nu)\.\Jc\.Y_k(\nu)\Big|_{0}^\infty=-2\.i\.\nu\.\sum_{k=0}^\infty Y_k^*(\nu)\.\Ps_k\.Y_k(\nu)
  \end{align*}
 with the left-hand side equal to $-Y_0^*(\nu)\.\Jc\.Y_0(\nu)$ by the limit point case hypothesis and 
 Theorem~\ref{T:A3}, i.e.,
  \begin{equation*}%\label{E:MR22}
   \sum_{k=0}^\infty Y_k^*(\nu)\.\Ps_k\.Y_k(\nu)=(2\.i\.\nu)^{-1}\.Y_0^*(\nu)\.\Jc\.Y_0(\nu).
  \end{equation*}
 From the Laurent series expansion in~\eqref{E:MR5} we obtain
  \begin{equation*}
   Y_0(\nu)=\al^*-\Jc\.\al^*\.\big(K_0+i\.\nu\.K_1+(i\.\nu)^2\.K_2+\dots\big),
  \end{equation*}
 and so 
  \begin{equation}\label{E:MR23}
   \sum_{k=0}^\infty Y_k^*(\nu)\.\Ps_k\.Y_k(\nu)=K_1+(i\.\nu)^2\.K_3+(i\.\nu)^4\.K_5+\dots=K_1+\O(\nu^2)
  \end{equation}
 for all $\nu$ small enough, i.e., the sum $\sum_{k=0}^m Y_k^*(\nu)\.\Ps_k\.Y_k(\nu)$ is uniformly bounded. Since
  \begin{align*}
   Y_k^+(\la_0)\coloneq\lim_{\nu\to0^+} Y_k(\nu)
    &=\hZ_k(\la_0)+\tZ_k(\la_0)\.K_0 \notag \\
    &\hspace{10mm} +\lim_{\nu\to0^+}\big[(i\.\nu)^{-1}\big(\tZ_k(\la_0+i\.\nu)-\tZ_k(\la_0)\big)\.K_{-1}\\
                           &\hspace*{40mm}+i\.\nu\.\tZ_k(\la_0+i\.\nu)\.K_1+\dots\big ]\\ 
    &=\hZ_k(\la_0)+\tZ_k(\la_0)\.K_0+\Big[\frac{\d}{\d\la} \tZ_k(\la)\Big]_{\la=\la_0}\.K_{-1}, %\label{E:MR24}
  \end{align*}
 it follows from the uniform boundedness in~\eqref{E:MR23} and the Moore--Osgood theorem that
  \begin{align*}
   \normP{Y^+(\la_0)}^2
    &=\lim_{m\to\infty}\.\sum_{k=0}^m Y_k^{+*}(\la_0)\.\Ps_k\.Y_k^+(\la_0)
     =\lim_{m\to\infty}\lim_{\nu\to0^+}\.\sum_{k=0}^m Y_k^*(\nu)\.\Ps_k\.Y_k(\nu)\\
%     =\lim_{\nu\to0}\sum_{k=0}^\infty Y_k^*\.\Ps_k\.Y_k
    &=\lim_{\nu\to0^+} \lim_{m\to\infty}\.\sum_{k=0}^m  Y_k^*(\nu)\.\Ps_k\.Y_k(\nu)
     \leq \lim_{\nu\to0^+} \lim_{m\to\infty} \big(K_1+\O(\nu^2)\big)=K_1.
  \end{align*}
 It means that the columns of $Y^+(\la_0)$ belong to $\ltp$ and the proof is complete.
\end{proof}

It remains to characterize the elements of the essential spectrum.

\begin{theorem}\label{T:MR3}
 Let Hypothesis~\ref{H:basic} be satisfied. We have $\la_0\in\sipc(\TLP)$ if and only if $M_+(\la)$ is not holomorphic 
 at $\la_0$, it holds
  \begin{equation}\label{E:MR26}
   L\coloneq\lim_{\nu\to0} \nu\.M_+(\la_0+i\.\nu)\neq0,
  \end{equation}
 and $M_+(\la)-i\.L\.(\la-\la_0)^{-1}$ is not holomorphic at $\la_0$. In this case, the columns of 
 $\tZ(\la_0)\.L$ belong to $\ltp$ and its nonzero columns are eigenfunctions of $\TLP$ related to $\la_0$.
\end{theorem}
\begin{proof}
 Let $M_+(\la)$ not be holomorphic at $\la_0$, equality~\eqref{E:MR26} hold, and $M_+(\la)-i\.L\.(\la-\la_0)^{-1}$ not 
 be holomorphic at $\la_0$. Obviously, $\la_0\not\in\rh(\TLP)$ by Theorem~\ref{T:MR1}, so $\la_0\in\si(\TLP)$. We 
 cannot have $\la_0\in\sip(\TLP)$ either, because otherwise $\la_0$ would be a~simple pole and $L=-i\.K_{-1}$ by 
 Theorem~\ref{T:MR2} and Lemma~\ref{L:MR2A}, which would imply
  \begin{equation*}
   M_+(\la)-i\.L\.(\la-\la_0)^{-1}=K_0+K_1\.(\la-\la_0)+K_2\.(\la-\la_0)^2+\dots,
  \end{equation*}
 i.e., $M_+(\la)-i\.L\.(\la-\la_0)^{-1}$ would be holomorphic at $\la_0$ contradicting our original assumption. Thus 
 $\la_0\in\sie(\TLP)$ and it remains to show the square summability of the columns of $\tZ(\la_0)\.L$, because in that 
 case the equality $\al\.\tZ_0(\la_0)\.L\.\xi=0$ for any $\xi\in\Cbb^n$ guarantees that the nonzero columns of 
 $\tZ_0(\la_0)\.L$ are eigenfunctions of $\TLP$ corresponding to $\la_0$, i.e., $\la_0\in\sipc(\TLP)$. For 
 $\la=\la_0+i\.\nu$ and any $m\in\oinftyZ$ we have
  \begin{equation*}
   \nu^2\.\sum_{k=0}^m \Xc_k^{+*}(\la)\.\Ps_k\.\Xc_k^{+}(\la)
    \leq \frac{\nu}{2i}\.\big(M_+(\la_0+i\.\nu)-M_+(\la_0-i\.\nu)\big).
  \end{equation*}
 From~\eqref{E:MR26} it follows also $\lim_{\nu\to0} \nu\.M(\la_0-i\.\nu)=L^*$, thereby
  \begin{align*}
   \sum_{k=0}^m \big(\tZ_k(\la_0)\.L\big)^*\.\Ps_k\.\big(\tZ_k(\la_0)\.L\big)
    &=\lim_{\nu\to0} 
     \nu^2\.\sum_{k=0}^m \Xc_k^{+*}(\la)\.\Ps_k\.\Xc_k^{+}(\la)\\
    &\hspace*{-15mm}\leq \frac{1}{2i} \lim_{\nu\to0} \big(\nu\.(M_+(\la_0+i\.\nu)-\nu\.M_+(\la_0-i\.\nu)\big)
%     =\frac{1}{2i} (L-L^*)
    =\im L.
  \end{align*}
 Similarly as above, this uniform boundedness implies
  \begin{equation*}
   \sum_{k=0}^\infty \big(\tZ_k(\la_0)\.L\big)^*\.\Ps_k\.\big(\tZ_k(\la_0)\.L\big)\leq \im L,
  \end{equation*}
 i.e., the columns of $\tZ(\la_0)\.L$ are square summable.
 
 On the other hand, let $\la_0\in\sipc(\TLP)$. From Lemma~\ref{L:MR2A} we know that 
 $L\coloneq \lim_{\nu\to0} \nu\.M(\la_0+i\.\nu)=i\.\big(\tau(\la_0)-\tau(\la^-)\big)$ exists. Let us assume that $L=0$. 
 If $z(\la_0)\in\Cbb(\oinftyZ)^{2n}$ is an eigenfunction and $\la=\la_0+i\.\nu$ for some $\nu\in\Rbb\stm\{0\}$, then 
 for $\Xc^+(\la_0+i\.\nu)=\big(\Xc^{[1]+}(\la_0+i\.\nu),\,\dots,\,\Xc^{[n]+}(\la_0+i\.\nu)\big)$ we have 
 $[z(\la_0)],[\Xc^{[1]+}(\la_0+i\.\nu)],\dots,[\Xc^{[n]+}(\la_0+i\.\nu)]\in\dom \Tmax$. Since~\eqref{E:Sla} is in the 
 limit point case, it follows
  \begin{equation*}
   \lim_{k\to\infty} \Xc_k^{[j]*+}(\la_0+i\.\nu)\.\Jc\.z_k(\la_0)=0 \qtext{for all $ j\in\{1,\dots,n\}$}
  \end{equation*}
 by Theorem~\ref{T:A3}, and so
  \begin{equation*}
   -i\.\nu\.\sum_{k=0}^\infty \Xc_k^{[j]*+}(\la_0+i\.\nu)\.\Ps_k\.z_k(\la_0)
    =-\Xc_0^{[j]*+}(\la_0+i\.\nu)\.\Jc\.z_0(\la_0)=\xi_1^{[j]},
  \end{equation*}
 where $\xi=(\xi_1,\dots,\xi_n)^\top\in\Cbb^n\stm\{0\}$ is such that $z(\la_0)=-\tZ(\la_0)\.\xi$ 
 by~\eqref{E:hz.tz.def}. Hence, by the Cauchy--Schwartz inequality, the Lagrange identity, and the 
 first part of~\eqref{E:MR15.i} we obtain
  \begin{align*}
   0\neq \xi^*\.\xi
    &=\abs{\nu}\times\abs[\Big]{\sum_{k=0}^\infty \xi^*\.\Xc_k^{*+}(\la_0+i\.\nu)\.\Ps_k\.z_k(\la_0)}\\
    &\leq \abs{\nu}\times\normP{\Xc_k^{+}(\la_0+i\.\nu)\.\xi}\times\normP{z(\la_0)}
     =\big(\abs{\nu}\times\xi^*\.\im M_+(\la_0+i\.\nu)\.\xi\big)^{1/2}\.\.\normP{z(\la_0)}.
  \end{align*}
 However, this leads to the contradiction, because the right-hand side tends to zero as $\nu\to0$ by the 
 assumption $L=0$. Therefore, equality~\eqref{E:MR26} holds and, consequently, $M_+(\la)$ is not holomorphic at 
 $\la_0$ as discussed before Lemma~\ref{L:MR2C}, see also Lemma~\ref{L:MR2A}. Finally, the function 
 $M_+(\la)-i\.L\.(\la-\la_0)^{-1}$ cannot be holomorphic at $\la_0$, because otherwise $\la_0$ would be a~simple pole 
 of $M_+(\la)$ implying $\la_0\in\sip(\TLP)$ by Theorem~\ref{T:MR2}.
\end{proof}

\begin{theorem}\label{T:MR4}
 Let Hypothesis~\ref{H:basic} be satisfied. We have $\la_0\in\sic(\TLP)$ if and only if $M_+(\la)$ is not holomorphic 
 at $\la_0$ and 
  \begin{equation*}
   \lim_{\nu\to0} \nu\.M_+(\la_0+i\.\nu)=0.
  \end{equation*}
\end{theorem}
\begin{proof}
 Let $\la_0\in\sic(\TLP)$. Then $\la_0\in\sie(\TLP)$ and $\la_0\not\in\sipc(\TLP)$ and $M_+(\la)$ is not holomorphic 
 at $\la_0$ by Theorem~\ref{T:MR1}. By Lemma~\ref{L:MR2A}, the limit $L\coloneq\lim_{\nu\to0} \nu\.M_+(\la_0+i\.\nu)$
 exists. Let us assume that $L\neq0$. Then $\tZ(\la_0)\.L\neq0$ and the columns of $\tZ(\la_0)\.L$ belong to $\ltp$ 
 by the first part of the proof of Theorem~\ref{T:MR3}, i.e., $\tZ(\la_0)\.L\.\xi\neq0$ is an eigenfunction of 
 $\TLP$ related to $\la_0$ for a suitable $\xi\in\Cbb^n\stm\{0\}$, because $\al\.\tZ_0(\la_0)\.L\.\xi=0$. But this
 is possible only when $\la_0\in\sip(\TLP)\cup\sipc(\TLP)$, which contradicts the original assumption, and so $L=0$.
 
 On the other hand, if $M_+(\la)$ is not holomorphic at $\la_0$, then $\la_0\not\in\rh(\TLP)$, i.e., 
 $\la_0\in\si(\TLP)$. If $\lim_{\nu\to0} \nu\.M_+(\la_0+i\.\nu)=0$, then $\la_0\not\in\sipc(\TLP)$ by 
 Theorem~\ref{T:MR3} and also $\la_0\not\in\sip(\TLP)$ by Theorem~\ref{T:MR2}, especially by~\eqref{E:MR6} and 
 Lemma~\ref{L:MR2A}. Therefore, $\la_0\in\sic(\TLP)$ by definition.
\end{proof}

This completes the proof of Theorem~\ref{T:intro.spectrum}.

\section{Applications}\label{S:application}

As a direct consequence of Theorems~\ref{T:MR1}--\ref{T:MR4} and Lemma~\ref{L:MR2C} we obtain the following 
conclusion. It is a partial extension of \cite{oD.rH03} for a~class of Sturm--Liouville difference equations, where 
a~necessary and sufficient condition for the discreteness and boundedness below of the spectrum of an associated 
difference operator was derived, see also \cite{dbH.rtL75,dbH.rtL78,dbH.rtL79,dtS91:PRSESA}. A more general result will 
be the main object of our subsequent research.
 
\begin{corollary} \label{C:MR2D}
 Let Hypothesis~\ref{H:basic} be satisfied. The linear relation $\TLP$ has a pure discrete spectrum if and only if 
 $M_+(\la)$ is meromorphic.
\end{corollary}

Furthermore, from the openness of $\rh(\TLP)$ and the definition of $\sip(\TLP)$ we obtain the first part of the 
following corollary, while its second part is a consequence of Lemmas~\ref{L:MR2A}--\ref{L:MR2C}, 
Theorems~\ref{T:MR2}--\ref{T:MR3}, and the monotone property of the function $\tau(\cdot)$.

\begin{corollary}\label{C:sets}
 Let Hypothesis~\ref{H:basic} be satisfied. The set $\rh(\TLP)\cup\sip(\TLP)$ is open and its complement 
 $\sie(\TLP)$ is closed. Furthermore, the set $\sip(\TLP)\cup\sipc(\TLP)$ is at most countable. 
\end{corollary}

All of the previous results were established for a fixed $\al\in\Ga$, which determines the initial values of the 
fundamental solutions $\hZ(\la)$ and $\tZ(\la)$ and of sequences in the domain of $\TLP$. How does its change affect 
the spectrum? This is investigated in the following statement, in which we emphasize the dependence of $\TLP$ on $\al$ 
as $\TLP(\al)$. It shows that the essential spectrum is invariant with respect to $\al$, while elements of the 
resolvent 
set and the discrete spectrum may alternate, compare with \cite{dbH.jkS82:QM,yS06}. It is based on the connection 
between two $M_+(\la)$-functions determined by different values of $\al$, which reads as
 \begin{equation}\label{E:MR27}  
  M_+(\la,\al)=\big[\al\.\Jc\.\hal^*+\al\.\hal^*\.M_+(\la,\hal)\big]\times
                       \big[\al\.\hal^*-\al\.\Jc\.\hal^*\.M_+(\la,\hal)\big]^{-1}
 \end{equation}
for any $\al,\hal\in\Ga$. Note that the proof of~\eqref{E:MR27} is similar as in~\cite[Theorem~2.13]{rSH.pZ14:JDEA}, 
because the limit point case guarantees that the columns of the corresponding Weyl solutions $\Xc^+(\la,M_+(\la,\al))$ 
and $\Xc^+(\la,M_+(\la,\hal))$ span the same $n$-dimensional space of square summable solutions of~\eqref{E:Sla} for 
any $\la\in\Cbb\stm\Rbb$.

\begin{theorem}\label{T:MR5}
 Let Hypothesis~\ref{H:basic} be satisfied. The sets $\rh(\TLP(\al))\cup\sip(\TLP(\al))$ and $\sie(\TLP(\al))$ are 
 invariant with respect to $\al$.
\end{theorem}
\begin{proof}
 Let $\hal\in\Ga$ be arbitrary. Then $\Cbb\stm\Rbb\subseteq\rh(\TLP(\hal))$. If $\la_0\in\rh(\TLP(\al))\cap\Rbb$, then 
 $M_+(\la,\al)$ is holomorphic for all $\la\in\Oc(\la_0)$ and the same is true for 
 $\al\.\Jc\.\hal^*+\al\.\hal^*\.M_+(\la,\hal)$ and $\al\.\hal^*-\al\.\Jc\.\hal^*\.M_+(\la,\hal)$. Consequently, the 
 latter function may have only isolated zeros on $\Oc(\la_0)$, which are poles of 
 $\big[\al\.\hal^*-\al\.\Jc\.\hal^*\.M_+(\la,\hal)\big]^{-1}$, so equality~\eqref{E:MR27} holds for all 
 $\la\in\Oc(\la_0)$ except for these poles by the Identity theorem for holomorphic function, i.e., we have either 
 $\la_0\in\rh(\TLP(\al))$ or $\la_0\in\sip(\TLP(\hal))$ by Lemma~\ref{L:MR2C} and Theorem~\ref{T:MR2}. If 
 $\la_0\in\sip(\TLP(\al))$, then $M_+(\la,\al)$ is holomorphic for all $\la\in\Oc^*(\la_0)$ and by similar arguments as 
 above we get $\la_0\in\rh(\TLP(\hal))\cup\sip(\TLP(\hal))$. This means that 
 $\rh(\TLP(\al))\cup\sip(\TLP(\al))\subseteq\rh(\TLP(\hal))\cup\sip(\TLP(\hal))$ and reversing the roles of the 
 matrices $\al$ and $\hal$ we get the opposite inclusion, which yields the independence of 
 $\rh(\TLP(\al))\cup\sip(\TLP(\al))$ on the choice of $\al$. Finally, the independence of $\sie(\TLP(\al))$ 
 is a straightforward consequence of the previous part.
\end{proof}

The Sturmian theory of discrete symplectic systems belongs to the most studied topics in this field and a 
comprehensive summary of the current state of this research can be found in~\cite[Chapters~4 and~5]{oD.jE.rSH19}. 
In the final part of this section, we provide an alternative approach to this theory on the interval $\oinftyZ$, 
compare with~\cite[Section~6.4]{oD.jE.rSH19} and~\cite[Subsection~7.3]{yS06}.

\begin{theorem}\label{T:MR6}
 Let Hypothesis~\ref{H:basic} be satisfied and $(a,b)\subseteq\rh(\TLP(\al))$. The for any $\hal\in\Ga$ we have
 $(a,b)\subseteq\rh(\TLP(\hal))$ except for at most $m\coloneq \rank\al\.\Jc\.\hal^*$ points belonging to 
 $\sip(\TLP(\hal))$.
\end{theorem}
\begin{proof}
 By Theorem~\ref{T:MR1} and Lemma~\ref{L:MR2C}, the function $M_+(\la)$ is holomorphic for all $\la\in(a,b)$ and we aim 
 to show that $M_+(\la,\hal)$ possesses at most $m$ simple poles in the interval $(a,b)$. From~\eqref{E:MR27} we can 
 see that only the zeros of the determinant of $\al\.\hal^*-\al\.\Jc\.\hal^*\.M_+(\la,\hal)$ may be singular points of 
 $M_+(\la,\hal)$. Let us put $\om\coloneq(\al\.\hal^*,\,\,\al\.\Jc\.\hal^*)$. Then $\om\in\Ga$, because 
 $(\hal^*,\,\,-\Jc\.\hal^*)$ is a unitary matrix by~\cite[Equality~(2.7)]{rSH.pZ14:JDEA}. Moreover, the equality 
 $\rank\al\.\Jc\.\hal^*=m$ implies the existence of the singular value decomposition of $\al\.\Jc\.\hal^*$, i.e., 
 $\al\.\Jc\.\hal^*=U\.\diag\{Q,0\}\.V$ with $n\times n$ unitary matrices $U,V$ and $Q\in\Cbb^{m\times m}$ such that
 $\rank Q=m$. Furthermore, also for
  \begin{equation*}
   \tom\coloneq (\tom_1,\,\,\tom_2) \qtextq{with} 
   \tom_1=\msmatrix{\tom_{11} & \tom_{12}\\ \tom_{21} & \tom_{22}}\coloneq U^*\.\al\.\hal^*\.V 
  \end{equation*}
 and
  \begin{equation*}
   \tom_2\coloneq U^*\.\al\.\Jc\.\hal^*\.V=\diag\{Q,0\}                
  \end{equation*}
 we have $\tom\in\Ga$ and, consequently, the blocks $\tom_{11}(\la)\in\Cbb^{m\times m}$, 
 $\tom_{12}(\la)\in\Cbb^{m\times (n-m)}$, $\tom_{21}(\la)\in\Cbb^{(n-m)\times m}$, and 
 $\tom_{22}\in\Cbb^{(n-m)\times(n-m)}$ satisfy
  \begin{gather*}
   \tom_{11}\.Q^*=Q\.\tom_{11}^*,\quad Q\.\tom_{21}^*=0, \quad 
   \tom_{11}\.\tom_{11}^*+\tom_{12}\.\tom_{12}^*+Q\.Q^*=I_m,\\
   \tom_{12}\.\tom_{22}^*=0, \qtextq{and} \tom_{22}\.\tom_{22}^*=I_{n-m}.
  \end{gather*}
 Therefore, 
  \begin{align*}
   \det\big(\al\.\hal^*-\al\.\Jc\.\hal^*\.M_+(\la,\hal)\big)
    &=\det\big(U\.\big(\tom_1-\tom_2\.K(\la)\big)\.V^*\big)\\
    &\hspace*{-10mm}
     =\det U\times \overline{\det V}\times \det\mmatrix{\om_{11}-Q\.K_{11}(\la) & -Q\.K_{12}(\la)\\ 0 & \tom_{22}}\\
    &\hspace*{-10mm}
     =\det U\times \overline{\det V}\times \det\tom_{22}\times \det(\om_{11}-Q\.K_{11}(\la))
  \end{align*}
 and its zeros are the same as for the determinant of the Hermitian matrix $\om_{11}\.Q^*-Q\.K_{11}(\la)\.Q^*$, where
 $K(\al)\coloneq V^*\.M_+(\al,\hal)\.V$ is an $n\times n$ Hermitian matrix with the left-upper block 
 $K_{11}(\la)\in\Cbb^{m\times m}$ and the right-upper block $K_{12}(\la)\in\Cbb^{m\times(n-m)}$. From 
 Lemma~\ref{L:MR1B}, it follows $K'(\la)>0$ for all $\la\in(a,b)$,  which means that 
 $\om_{11}\.Q^*-Q\.K_{11}(\la)\.Q^*$ is decreasing on $(a,b)$ and, consequently, its eigenvalues are 
 also decreasing. Thus, $\det\big(\al\.\hal^*-\al\.\Jc\.\hal^*\.M_+(\la,\hal)\big)$ may have at most $m$ zeros as 
 claimed above.
\end{proof}

In particular, if $\al\.\Jc\.\hal^*=0$, then $\al\.\hal^*$ is invertible by the previous proof and 
$M_+(\la,\al)=\al\.\hal^*\.M_+(\la,\hal)\.\big(\al\.\hal^*\big)^{-1}$. Therefore, the corresponding linear relations 
have the same resolvent.

\begin{corollary}\label{C:MR7}
 Let Hypothesis~\ref{H:basic} be satisfied and $\hal\in\Ga$ be such that $\al\.\Jc\.\hal^*=0$. Then 
 $\rh(\TLP(\al))=\rh(\TLP(\hal))$.
\end{corollary}

From Theorem~\ref{T:MR6} we obtain also the invariance of the boundedness of the spectrum.

\begin{corollary}\label{C:MR10}
 Let Hypothesis~\ref{H:basic} be satisfied and the set $\si(\TLP(\al))$ be bounded from below (or from above) with 
 $\mu_{\min}\coloneq\inf\{\la\in\si(\TLP(\al))\}$ (or $\mu_{\max}\coloneq\sup\{\la\in\si(\TLP(\al))\}$). Then the 
 spectrum of $\TLP(\hal)$ is also bounded from below (or above) for any $\hal\in\Ga$ and, in addition, the intersection
 $\sip(\TLP(\hal))\cap(-\infty,\mu_{\min})$ (or the intersection $\sip(\TLP(\hal))\cap(\mu_{\max},\infty)$) contains at 
 most $m\coloneq\rank \al\.\Jc\.\hal^*$ points.
\end{corollary}

We already know that the set $\rh(\TLP(\al))\cup\sip(\TLP(\al))$ is invariant with respect to $\al$. But what can we 
say about the discrete spectrum itself?

\begin{theorem}\label{T:MR8}
 Let Hypothesis~\ref{H:basic} be satisfied and $\la_1<\la_2<\dots<\la_{n+1}$ be $n+1$ consecutive points of 
 $\sip(\TLP(\al))$. If $(\la_1,\la_{n+1})\cap\sie(\TLP(\al))=\emptyset$, then 
 $[\la_1,\la_{n+1}]\cap\sip(\TLP(\hal))\neq\emptyset$ for any $\hal\in\Ga$, i.e., the interval $[\la_1,\la_{n+1}]$ 
 contains at least one eigenvalue of the linear relation $\TLP(\hal)$.
\end{theorem}
\begin{proof}
 Let $(\la_1,\la_{n+1})\cap\sie(\TLP(\al))=\emptyset$ and suppose that it also holds
 $[\la_1,\la_{n+1}]\cap\sip(\TLP(\hal))=\emptyset$ for some $\hal\in\Ga$. Then 
 $[\la_1,\la_{n+1}]\cap\si(\TLP(\hal))=\emptyset$ by Theorem~\ref{T:MR5}, i.e., 
 $[\la_1,\la_{n+1}]\subseteq\rh(\TLP(\hal))$ and from the openness of this set it follows that
 $(\la_1-\eps,\la_{n+1}+\eps)\subseteq\rh(\TLP(\hal))$ for some $\eps>0$. However, by Theorem~\ref{T:MR6} there is at 
 most $m\leq n$ points of $\sip(\TLP(\al))$ in $(\la_1-\eps,\la_{n+1}+\eps)$, which contradicts the fact 
 $\la_1,\dots,\la_{n+1}\in\sip(\TLP(\al))$.
\end{proof}

Upon combining Theorems~\ref{T:MR6} and~\ref{T:MR8} we get the following characterization of the case of a~pure 
discrete spectrum, compare with~\cite[Section~6.2.2]{oD.jE.rSH19}. Let us emphasize that collocation ``interlacing 
property'' and the word ``among'' are here used in the inclusive sense, which means that they do not exclude the left 
and right endpoints of the interval under consideration.

\begin{corollary}\label{C:MR9}
 Let Hypothesis~\ref{H:basic} be satisfied and $\TLP(\al)$ have pure discrete spectrum. Then the same is true for 
 $\TLP(\hal)$ for any $\hal\in\Ga$ and these spectra possesses the interlacing property, i.e., among any $n+1$ 
 consecutive points of $\sip(\TLP(\al))$ there is at least one element of $\sip(\TLP(\hal))$.
\end{corollary}

Finally, we focus on the scalar case $n=1$. In this case, every $1\times2$ matrix $\al\in\Ga$ determining the initial 
value for the sequences in $\dom\TLP$ can be written as
 \begin{equation}\label{E:alpha.n=1}
  \al=(\sin\al_0,\,\,\cos\al_0)
 \end{equation}
for some $\al_0\in[0,2\pi)$. If we take also $\hal=(\sin\hal_0,\,\,\cos\hal_0)$, then
 \begin{equation*}
  \al\.\hal^*=\cos(\al_0-\hal_0) \qtextq{and} \al\.\Jc\.\hal^*=\sin(\al_0-\hal_0),
 \end{equation*}
so, by~\eqref{E:MR27}, we have
 \begin{equation*}
  M_+(\la,\al)=M_+(\la,\hal) \qtext{whenever $\sin(\al_0-\hal_0)=0$,}
 \end{equation*}
and otherwise
 \begin{equation*}
  M_+(\la,\al)=[1+\cotan(\al_0-\hal_0)\.M_+(\la,\hal)]\times[\cotan(\al_0-\hal_0)-M_+(\la,\hal)]^{-1}.
 \end{equation*}
This together with the inequality $M_+'(\la)>0$ from Lemma~\ref{L:MR1B} shows that $M_+(\la,\al)$ has precisely one 
simple pole between two consecutive poles of $M_+(\la,\hal)$. Therefore, Theorems~\ref{T:MR6} and~\ref{T:MR8} can be 
stated as follows.
 
\begin{corollary}\label{C:MR11}
 Let Hypothesis~\ref{H:basic} be satisfied with $n=1$ and $\al$ given as in~\eqref{E:alpha.n=1} for some 
 $\al_0\in[0,2\pi)$. Furthermore, let $\la_1<\la_2$ be two consecutive points of $\sip(\TLP(\al))$ and 
 $\hal=(\sin\hal_0,\,\,\cos\hal_0)$ with $\hal_0\in[0,2\pi)$ be arbitrary.
  \begin{enumerate}[leftmargin=8mm,topsep=1mm,label={{\normalfont{(\roman*)}}}]
   \item If $\sin(\al_0-\hal_0)=0$, then $\rh(\TLP(\al))=\rh(\TLP(\hal))$ and the spectra $\si(\TLP(\al))$ and 
         $\si(\TLP(\hal))$ possess the same structure.

   \item On the other hand, if $\sin(\al_0-\hal_0)\neq0$ and $(\la_1,\la_2)\subseteq\rh(\TLP(\al))$ or equivalently
         $(\la_1,\la_2)\cap\sie(\TLP(\al))=\emptyset$, then $(\la_1,\la_2)\cap\sip(\TLP(\hal))$ is a singleton. 
         Especially, if $\TLP(\al)$ has a~pure discrete spectrum, then $\TLP(\hal)$ does the same and between 
         two consecutive points of $\sip(\TLP(\al))$ always lies precisely one element of $\sip(\TLP(\hal))$.
  \end{enumerate}
\end{corollary}

%%%%%%%%%%%%%%%%%%%%%%%%%%%%%%%%%%%%%%%%%%%% SECTION %%%%%%%%%%%%%%%%%%%%%%%%%%%%%%%%%%%%%%%%%%%%%%%%%%%%%%%%%%%%%%%%%

\section*{Acknowledgments}
\addcontentsline{toc}{section}{Acknowledgments}
The research was supported by the Czech Science Foundation under Grant GA19-01246S. The author is grateful to anonymous 
referees for detailed reading of the manu\-script and their comments and suggestions.

%%%%%%%%%%%%%%%%%%%%%%%%%%%%%%%%%%%%%%%%%%%% SECTION %%%%%%%%%%%%%%%%%%%%%%%%%%%%%%%%%%%%%%%%%%%%%%%%%%%%%%%%%%%%%%%%%

% \bibliographystyle{my_bibliography_paper_style}
% \bibliography{bibliotheca}

\begin{thebibliography}{10}
\addcontentsline{toc}{section}{References}

\bibitem{joA.doA.foN13}
J.~O. Agure, D.~O. Ambogo, and F.~O. Nyamwala, \textit{Deficiency indices and
  spectrum of fourth order difference equations with unbounded coefficients},
  Math. Nachr.~\textbf{286} (2013), no.~4, 323--339.

\bibitem{hB13}
H.~Behncke, \textit{Spectral theory of Hamiltonian difference systems with
  almost constant coefficients}, J. Difference Equ. Appl.~\textbf{19} (2013),
  no.~1, 1--12.

\bibitem{jB.sH.hsvdS20}
J.~Behrndt, S.~Hassi, and H.~S.~V. de~Snoo, \textit{Boundary Value Problems,
  Weyl Functions, and Differential Operators}, Monographs in Mathematics,
  Vol.~108, Birkhäuser, Cham, 2020. ISBN 978-3-030-36713-8; 978-3-030-36714-5.

\bibitem{rB.jE.aK.gT14}
R.~Brunnhuber, J.~Eckhardt, A.~Kostenko, and G.~Teschl, \textit{Singular
  Weyl--Titchmarsh--Kodaira theory for one-dimensional Dirac operators},
  Monatsh. Math.~\textbf{174} (2014), no.~4, 515--547.

\bibitem{jC.wnE68}
J.~Chaudhuri and W.~N. Everitt, \textit{On the spectrum of ordinary second
  order differential operators}, Proc. R. Soc. Edinb., Sect. A (1967-1968),
  95--119 (1969).

\bibitem{slC.pZ10}
S.~L. Clark and P.~Zem{\'{a}}nek, \textit{On a Weyl--Titchmarsh theory for
  discrete symplectic systems on a half line}, Appl. Math. Comput.~\textbf{217}
  (2010), no.~7, 2952--2976.

\bibitem{slC.pZ15}
S.~L. Clark and P.~Zem{\'{a}}nek, \textit{On discrete symplectic systems:
  Associated maximal and minimal linear relations and nonhomogeneous problems},
  J. Math. Anal. Appl.~\textbf{421} (2015), no.~1, 779--805.

\bibitem{oD.jE.rSH19}
O.~Do{\v{s}}l{\'{y}}, J.~Elyseeva, and R.~{\v{S}}imon~Hilscher,
  \textit{Symplectic Difference Systems: Oscillation and Spectral Theory},
  Pathways in Mathematics, Birkh{\"{a}}user/Springer, Cham, 2019. ISBN
  978-3-030-19372-0; 978-3-030-19373-7.

\bibitem{oD.rH03}
O.~Do{\v{s}}l{\'y} and R.~Hilscher, \textit{A class of Sturm--Liouville
  difference equations: (Non)oscillation constants and property BD}, in
  ``Advances in Difference Equations, IV'' (R.~P.~Agarwal, M.~Bohner, and
  D.~O'Regan, editors), Comput. Math. Appl.~\textbf{45} (2003), no.~6-9,
  961--981.

\bibitem{fG.erT00}
F.~Gesztesy and E.~R. Tsekanovski{\u\i}, \textit{On matrix-valued Herglotz
  functions}, Math. Nachr.~\textbf{218} (2000), 61--138.

\bibitem{dbH.rtL75}
D.~B. Hinton and R.~T. Lewis, \textit{Discrete spectra criteria for singular
  differential operators with middle terms}, Math. Proc. Cambridge Philos.
  Soc.~\textbf{77} (1975), 337--347.

\bibitem{dbH.rtL78}
D.~B. Hinton and R.~T. Lewis, \textit{Spectral analysis of second order
  difference equations}, J. Math. Anal. Appl.~\textbf{63} (1978), no.~2,
  421--438.

\bibitem{dbH.rtL79}
D.~B. Hinton and R.~T. Lewis, \textit{Singular differential operators with
  spectra discrete and bounded below}, Proc. Roy. Soc. Edinburgh Sect.
  A~\textbf{84} (1979), no.~1-2, 117--134.

\bibitem{dbH.jkS82:QM}
D.~B. Hinton and J.~K. Shaw, \textit{On the spectrum of a singular Hamiltonian
  system}, Quaestiones Math.~\textbf{5} (1982/83), no.~1, 29--81.

\bibitem{gaM.aS.mT19}
G.~A. Monteiro, A.~Slav{\'{i}}k, and M.~Tvrd{\'{y}},
  \textit{Kurzweil--Stieltjes integral: Theory and Applications}, Series in
  Real Analysis, Vol.~15, World Scientific Publishing, Hackensack, 2019. ISBN
  978-981-4641-77-7.

\bibitem{foN17}
F.~O. Nyamwala, \textit{Essential and continuous spectrum of symmetric
  difference equations}, Math. Nachr.~\textbf{290} (2017), no.~17-18,
  2977--2991.

\bibitem{wR76}
W.~Rudin, \textit{Principles of Mathematical Analysis}, third edition,
  McGraw--Hill Book, New York, 1976. International Series in Pure and Applied
  Mathematics.

\bibitem{yS06}
Y.~Shi, \textit{Weyl--Titchmarsh theory for a class of discrete linear
  Hamiltonian systems}, Linear Algebra Appl.~\textbf{416} (2006), no.~2-3,
  452--519.

\bibitem{rSH.pZ14:ICDEA}
R.~{\v{S}}imon~Hilscher and P.~Zem{\'{a}}nek, \textit{Generalized Lagrange
  identity for discrete symplectic systems and applications in Weyl--Titchmarsh
  theory}, in ``Theory and Applications of Difference Equations and Discrete
  Dynamical Systems'', Proceedings of the 19th International Conference on
  Difference Equations and Applications (Muscat, 2013), Z.~AlSharawi,
  J.~Cushing, and S.~Elaydi (editors), Springer Proceedings in Mathematics \&
  Statistics, Vol. 102, pp. 187--202, Springer, Berlin, 2014.

\bibitem{rSH.pZ14:AMC}
R.~{\v{S}}imon~Hilscher and P.~Zem{\'{a}}nek, \textit{Limit point and limit
  circle classification for symplectic systems on time scales}, Appl. Math.
  Comput.~\textbf{233} (2014), 623--646.

\bibitem{rSH.pZ14:JDEA}
R.~{\v{S}}imon~Hilscher and P.~Zem{\'{a}}nek, \textit{Weyl--Titchmarsh theory
  for discrete symplectic systems with general linear dependence on spectral
  parameter}, J. Difference Equ. Appl.~\textbf{20} (2014), no.~1, 84--117.

\bibitem{dtS91:PRSESA}
D.~T. Smith, \textit{On the spectral analysis of selfadjoint operators
  generated by second order difference equations}, Proc. Roy. Soc. Edinburgh
  Sect. A~\textbf{118} (1991), no.~1-2, 139--151.

\bibitem{hS.qK.yS16}
H.~Sun, Q.~Kong, and Y.~Shi, \textit{Essential spectrum of singular discrete
  linear Hamiltonian systems}, Math. Nachr.~\textbf{289} (2016), no.~2-3,
  343--359.

\bibitem{hS.yS14}
H.~Sun and Y.~Shi, \textit{Spectral properties of singular discrete linear
  Hamiltonian systems}, J. Difference Equ. Appl.~\textbf{20} (2014), no.~3,
  379--405.

\bibitem{hS.yS15}
H.~Sun and Y.~Shi, \textit{On essential spectra of singular linear Hamiltonian
  systems}, Linear Algebra Appl.~\textbf{469} (2015), 204--229.

\bibitem{jdW51}
J.~D. Weston, \textit{Inequalities for Riemann--Stieltjes integrals}, Math.
  Z.~\textbf{54} (1951), 272--274.

\bibitem{pZ20}
P.~Zem{\'{a}}nek, \textit{Linear operators associated with differential and
  difference systems: What is different?}, in ``Progress on Difference
  Equations and Discrete Dynamical Systems'', Proceedings of the International
  Conference on Differential \& Difference Equations and Applications 2019
  (London, 2019), S.~Baigent, M.~Bohner, and S.~Elaydi (editors), Springer
  Proceedings in Mathematics \& Statistics, Vol.~341, pp. 435--448, Springer,
  Berlin, 2020.

\bibitem{pZ21}
P.~Zem{\'{a}}nek, \textit{Eigenfunctions expansion for discrete symplectic
  systems with general linear dependence on spectral parameter}, J. Math. Anal.
  Appl.~\textbf{499} (2021), no.~2, Article no.~125054, 1--37 pp. (electronic).

\bibitem{pZ?:MN}
P.~Zem{\'{a}}nek, \textit{Non-limit-circle and limit-point criteria for
  symplectic dynamic systems on time scales}, Math. Nachr., to appear.

\bibitem{pZ.slC16}
P.~Zem{\'{a}}nek and S.~L. Clark, \textit{Characterization of self-adjoint
  extensions for discrete symplectic systems}, J. Math. Anal.
  Appl.~\textbf{440} (2016), no.~1, 323--350.

\bibitem{pZ.slC:SAE2}
P.~Zem{\'{a}}nek and S.~L. Clark, \textit{Discrete symplectic systems, boundary
  triplets, and self-adjoint extensions}, submitted, 2020.

\end{thebibliography}

\end{document}